\documentclass[11pt,letterpaper]{amsart}

\usepackage[T1]{fontenc}
\usepackage{lmodern}
\usepackage[expansion=false]{microtype}
\usepackage{mathtools}
\usepackage{amssymb}
\usepackage{mathrsfs}
\usepackage{enumitem}
\usepackage{aliascnt}
\usepackage{needspace}
\usepackage[margin=1.08in]{geometry}
\usepackage[colorlinks=true,linkcolor=blue,citecolor=red,urlcolor=cyan,pagebackref=true]{hyperref}
\usepackage[capitalise,noabbrev]{cleveref}

\numberwithin{equation}{section}
\newtheorem{thm}{Theorem}[section]
\newtheorem*{problem}{Problem}
\newaliascnt{prop}{thm}
\newtheorem{prop}[prop]{Proposition}
\aliascntresetthe{prop}
\newaliascnt{lem}{thm}
\newtheorem{lem}[lem]{Lemma}
\aliascntresetthe{lem}
\newaliascnt{cor}{thm}
\newtheorem{cor}[cor]{Corollary}
\aliascntresetthe{cor}
\theoremstyle{definition}
\newaliascnt{defn}{thm}

\aliascntresetthe{defn}
\theoremstyle{remark}
\newaliascnt{remark}{thm}
\newtheorem{remark}[remark]{Remark}
\aliascntresetthe{remark}
\crefname{thm}{Theorem}{Theorems}
\Crefname{thm}{Theorem}{Theorems}
\crefname{prop}{Proposition}{Propositions}
\Crefname{prop}{Proposition}{Propositions}
\crefname{lem}{Lemma}{Lemmas}
\Crefname{lem}{Lemma}{Lemmas}
\crefname{cor}{Corollary}{Corollaries}
\Crefname{cor}{Corollary}{Corollaries}
\crefname{defn}{Definition}{Definitions}
\Crefname{defn}{Definition}{Definitions}
\crefname{remark}{Remark}{Remarks}
\Crefname{remark}{Remark}{Remarks}

\AddToHook{env/thm/before}{\Needspace{4\baselineskip}}
\AddToHook{env/prop/before}{\Needspace{4\baselineskip}}
\AddToHook{env/lem/before}{\Needspace{4\baselineskip}}
\AddToHook{env/cor/before}{\Needspace{4\baselineskip}}
\AddToHook{env/defn/before}{\Needspace{4\baselineskip}}
\AddToHook{env/remark/before}{\Needspace{4\baselineskip}}

\newcommand{\C}{\mathbb C}
\newcommand{\R}{\mathbb R}
\newcommand{\BB}{\mathbb B}
\newcommand{\DD}{\mathbb D}
\newcommand{\HH}{\mathbb H}
\newcommand{\Om}{\Omega}
\newcommand{\id}{\operatorname{id}}
\newcommand{\proofstep}[2]{\textbf{Step #1. #2.}\enspace}
\newcommand{\rank}{\operatorname{rank}}
\newcommand{\ran}{\operatorname{ran}}
\newcommand{\tr}{\operatorname{tr}}
\newcommand{\Aut}{\operatorname{Aut}}
\newcommand{\Rea}{\operatorname{Re}}
\newcommand{\Ima}{\operatorname{Im}}
\newcommand{\dist}{\operatorname{dist}}
\newcommand{\norm}[1]{\lVert#1\rVert}
\newcommand{\abs}[1]{\lvert#1\rvert}
\newcommand{\ip}[2]{\langle#1,#2\rangle}

\hypersetup{
  pdftitle={Normal Squeezing and Strict Pseudoconvexity},
  pdfsubject={Normal squeezing, local boundary geometry, and Bergman kernel growth},
  pdfcreator={},
  pdfproducer={},
  pdfkeywords={squeezing function, strict pseudoconvexity, Bergman kernel, holomorphic semigroup}
}

\begin{document}
\raggedbottom

\title[Normal squeezing and strict pseudoconvexity]
{Normal Limits of the Squeezing Function and Strict Pseudoconvexity}

\author[X. Qin]{Xiangsen Qin}
\address{Xiangsen Qin: Chern Institute of Mathematics and LPMC, Nankai University \\
Tianjin 300071, China}
\email{qinxiangsen@nankai.edu.cn}

\begin{abstract}
Let $p$ be a boundary point of a bounded pseudoconvex domain in
$\C^n$, $n\ge2$, and suppose the squeezing function tends to $1$
along the inward normal at $p$. We show that $p$ is strictly
pseudoconvex if the boundary is $C^{2,\alpha}$ near $p$ for some
$0<\alpha\le1$, without assuming finite type or convexity.
The proof begins under local $C^2$ regularity near $p$.
A comparison of Bergman kernels rules out nonconstant holomorphic
curve germs through $p$ contained in the boundary. With these
curve germs excluded, rescaled normal translations give rise to
a holomorphic semigroup. Classifying this semigroup identifies the normalized
kernel limit at selected scales and yields the logarithmic growth
exponent $n+1$ for the diagonal Bergman kernel along the inward normal.
Under the $C^{2,\alpha}$ assumption, Levi degeneracy at $p$ would
give a kernel upper bound with exponent strictly below $n+1$.
This contradiction proves $p$ is strictly pseudoconvex.
\end{abstract}


\maketitle
{\small\setlength{\baselineskip}{11pt}\tableofcontents}
\clearpage
\section{Introduction}\label{sec:intro}

Deng, Guan and Zhang introduced the squeezing function for bounded domains in \(\C^n\)
\cite[Definition~1.1]{DGZbasic}. It measures how closely a pointed domain
resembles the unit ball under holomorphic embeddings. For a
bounded domain $D\subset\C^n$ and $z\in D$, set
\begin{equation}\label{eq:squeezing}
 s_D(z)=\sup\bigl\{r\in(0,1):
 r\BB^n\subset f(D),\ f:D\hookrightarrow\BB^n
 \text{ holomorphic},\ f(z)=0\bigr\}.
\end{equation}
The same definition applies to any domain admitting a
holomorphic embedding into the ball. The function is biholomorphically
invariant. Deng, Guan and Zhang proved that it tends to one at the boundary
of a bounded domain with $C^2$ strictly pseudoconvex boundary
\cite[Theorem~1.3]{DGZ}.
Forn\ae ss--Wold asked whether this implication can be reversed.

\begin{problem}[Forn\ae ss--Wold {\cite[Problem~4.1, p.~85]{FW}}]
For a bounded pseudoconvex domain $\Om\subset\C^n$ with
$C^\infty$ boundary, does the condition
\[
 \lim_{\substack{z\in\Om\\\dist(z,\partial\Om)\to0}}
 s_\Om(z)=1
\]
force $\Om$ to be strictly pseudoconvex?
\end{problem}

Joo and Kim \cite[Question~1]{JK} record the local version,
attributed to Forn\ae ss: if squeezing tends to one at a specified
boundary point, must the Levi form there be positive definite?
Forn\ae ss and Wold also ask how much regularity is needed for the
global converse \cite[Problem~4.2, p.~85]{FW}. Their bounded convex
$C^2$ example shows that the converse fails at that regularity
\cite[Theorem~1.1]{FW}. With a local H\"older condition on the second
derivatives, we prove a converse using only the squeezing limit
along the inward normal.

\begin{thm}\label{thm:main}
Let $\Om\Subset\C^n$, $n\ge2$, be a bounded pseudoconvex domain.
Let $p\in\partial\Om$, and suppose that the boundary is
$C^{2,\alpha}$ in a neighborhood of $p$ for
some $0<\alpha\le1$. Let $\nu_p$ be the outward unit normal at $p$. If
\begin{equation}\label{eq:normal-hypothesis}
 \lim_{\delta\downarrow0}s_\Om(p-\delta\nu_p)=1,
\end{equation}
then $p$ is strictly pseudoconvex. No boundary regularity
away from this neighborhood is assumed.
\end{thm}

The limit in \eqref{eq:normal-hypothesis} is taken over all $\delta\downarrow0$,
not just a sequence. Both the full limit at $p$ and the global boundary
limit imply it. The local H\"older assumption can be replaced by a
pointwise bound of order $2+\alpha$ on the second-order
Taylor remainder of one local defining function
(\Cref{cor:pointwise-taylor}).

\begin{cor}\label{cor:local-question}
Let $\Om\Subset\C^n$, $n\ge2$, be a bounded pseudoconvex
domain, let $p\in\partial\Om$, and suppose that
$\partial\Om$ is $C^{2,\alpha}$ near $p$ for some
$0<\alpha\le1$. If
\[
 \lim_{\Om\ni z\to p}s_\Om(z)=1,
\]
then $p$ is strictly pseudoconvex. In particular, a bounded
pseudoconvex domain $\Om\Subset\C^n$, $n\ge2$, with
$C^{2,\alpha}$ boundary for some $0<\alpha\le1$ is strictly
pseudoconvex if $s_\Om(z)\to1$ as
$\dist(z,\partial\Om)\to0$, with $z\in\Om$.
\end{cor}

The regularity assumptions in \Cref{thm:main,cor:local-question}
cannot in general be weakened to $C^2$, as shown by the example
of Forn\ae ss and Wold \cite[Theorem~1.1]{FW}.

\Needspace{9\baselineskip}
At local $C^2$ regularity, the same normal squeezing limit still rules
out nonconstant holomorphic curve germs through $p$ contained in
$\partial\Om$ (\Cref{thm:no-curves}). It also determines the logarithmic
growth of the diagonal Bergman kernel $K_\Om(z)=K_\Om(z,z)$.

\begin{thm}\label{thm:normal-log-growth}
Let $\Om\Subset\C^n$, $n\ge2$, be a bounded pseudoconvex
domain whose boundary is $C^2$ near $p\in\partial\Om$ and let $\nu_p$ be the outward unit normal. If
\eqref{eq:normal-hypothesis} holds, then
\begin{equation}\label{eq:normal-log-growth}
 \lim_{\delta\downarrow0}
 \frac{\log K_\Om(p-\delta\nu_p)}{\log(1/\delta)}=n+1.
\end{equation}
\end{thm}

The logarithmic limit does not assert that
$\delta^{n+1}K_\Om(p-\delta\nu_p)$ tends to a positive finite number.

\Cref{cor:normal-sampling,cor:nontangential-curve} give two ways
to obtain the full normal limit from more restricted data. The
first uses a nontangential sequence with a quantitative condition
on adjacent depths and squeezing defects. That condition holds, for
example, for normal samples $p-2^{-j}\nu_p$ with squeezing tending
to one. The second uses a continuous nontangential curve,
such as a $C^1$ curve with an inward transverse derivative,
along which squeezing tends to one.
In both criteria, the remainder bound \eqref{eq:pointwise-taylor}
at $p$ can replace the local $C^{2,\alpha}$ assumption
(\Cref{cor:pointwise-taylor}).

\Cref{thm:main} also applies to biholomorphic realizations.
If $D\Subset\C^n$, $n\ge2$, is strictly pseudoconvex with
$C^2$ boundary and $\Om\Subset\C^n$ is biholomorphic to
$D$, then every $p\in\partial\Om$ at which the boundary is
locally $C^{2,\alpha}$ for some $0<\alpha\le1$ is strictly
pseudoconvex (\Cref{cor:strongly-pseudoconvex-realization}).
This consequence of \Cref{thm:main} requires no boundary
extension of the biholomorphism. The classical unbounded ball image
in \Cref{prop:unbounded-ball-image} has real-analytic boundary
with Levi-degenerate points, so boundedness cannot be dropped
for general targets.

\subsection{Related results}\label{subsec:earlier-results}

A smooth boundary point $p$ with defining function $\rho$ is of
finite D'Angelo type \cite{DA82} if
\[
 \sup_\gamma\frac{\operatorname{ord}_0(\rho\circ\gamma)}
 {\operatorname{ord}_0(\gamma-p)}<\infty,
\]
where $\operatorname{ord}_0$ denotes vanishing order at $0$
(the minimum of the component orders for a vector-valued map),
and $\gamma$ ranges over nonconstant holomorphic curve germs through $p$;
see \cite[Definitions~6.1--6.2, p.~213]{DK99}.
For the notion of weighted nontangential convergence,
see \cite[Definition~3.4]{NN20}.

The finite-type case has been studied under several approach
conditions. At a finite D'Angelo type point of a smoothly
bounded pseudoconvex domain in $\C^2$, Joo and Kim
\cite[Theorem~1.2]{JK} deduce strict pseudoconvexity from
the full squeezing limit at that point. At an $h$-extendible
point of a smoothly bounded pseudoconvex domain, Nikolov
\cite[Theorem~1]{Nik} deduces strict pseudoconvexity from
squeezing tending to one along a single nontangential sequence. Mahajan--Verma
\cite[Theorem~1.2 and the following paragraph]{MV} require
smoothness and $h$-extendibility only near the point in a bounded
domain and observe that a sequence along the inward normal suffices.
Ninh--Nguyen \cite[Theorem~1.2]{NN20} obtain a criterion for
weighted nontangential sequences in smoothly bounded pseudoconvex
domains.
Here $h$-extendibility means finite type together with equality
of the Catlin multitype, the weight invariant of adapted
coordinates \cite{Cat84}, and the D'Angelo multitype, the tuple of D'Angelo $q$-types
in reverse dimension order; see \cite[Definition~6.2; Section~7,
(24)--(25)]{DK99}. It includes
pseudoconvex finite-type points in $\C^2$, pseudoconvex
finite-type points of Levi corank at most one, and convex
finite-type points. These criteria impose no sampling
condition of the kind required in \Cref{cor:normal-sampling},
but their geometric assumptions and approach conditions differ
from those used here.

Bharali \cite[Theorem~1.3]{Bha} uses a single paraboloidal
sequence approaching a $C^\infty$ finite-type point of a
bounded domain in $\C^2$, with pseudoconvexity required only
near the point. Along that sequence, the norm of the real
tangential projection of the displacement is bounded by a
constant times the square root of the boundary distance.
Guermazi \cite{Gue25} studies finite-type domains with a locally
diagonalizable Levi form. Along tangential sequences approaching weakly pseudoconvex
points, the squeezing function can tend to one, as in
Ninh--Nguyen--Nguyen \cite[Example~1.2, pp.~832--833]{NNN}. Ninh \cite[Theorems~1.2 and~1.3]{Ninh26} studies the localized
invariant $s_{\Om\cap U}$ along specified tangential sequences
at finite-type points. These results do not provide a criterion for
the squeezing function of the full domain $\Om$ along the normal
without finite-type assumptions.

For convex domains, Zimmer's smooth gap theorem \cite[Theorem~1.2]{Zgap}
and its $C^{2,\alpha}$ version \cite[Theorem~1.7]{Zim} give
uniform quantitative characterizations. In the latter result,
a lower bound $1-\varepsilon$ for squeezing outside a
compact set suffices; $\varepsilon$ depends only on the
dimension and $\alpha$. The generalized-squeezing
criterion of Biswas--Chatterjee--Sarkar
\cite[Corollary~1.7]{BCS26} also assumes a bounded convex
$C^{2,\alpha}$ domain and uses a bounded balanced strongly
convex $C^3$ reference domain. \Cref{thm:main} concerns
one point of a general pseudoconvex domain and does not give
a domain-independent quantitative gap.

Nikolov--Verma \cite[Proposition~6]{NV} also obtain a local
converse at smooth locally convexifiable boundary points without
boundary analytic discs, assuming the full boundary squeezing
limit. Their result does not require finite type.

For biholomorphic realizations, Pinchuk \cite[Theorem~3.1]{Pin81}
proved that a bounded convex domain with $C^{2,\alpha}$
boundary, biholomorphic to a bounded $C^2$
strictly pseudoconvex domain is strictly pseudoconvex.
Bell \cite{Bell81} proved
invariance of strict pseudoconvexity for bounded $C^\infty$
domains. Zimmer \cite[Theorem~1.1 and Remark~1.2]{Zim}
showed that a convex $C^{2,\alpha}$ domain, $0<\alpha\le1$,
biholomorphic to a bounded $C^2$ strictly pseudoconvex domain is strictly
pseudoconvex; the convex target need not be bounded.
Hu \cite{Hu96} proved that a bounded
weakly but not strictly pseudoconvex domain with $C^k$ boundary
cannot be biholomorphic to a bounded strictly pseudoconvex domain
with $C^l$ boundary when $k>16n+15$ and $l>6n+8$.
\Cref{cor:strongly-pseudoconvex-realization} allows bounded
nonconvex targets with $C^{2,\alpha}$ boundary regularity assumed only near
the point under consideration. It does not address the $C^2$ ball-image
question in \cite[Problem~4.5]{FW}.

For a bounded pseudoconvex domain, H\"ormander
\cite[Theorem~3.5.1, p.~145]{Hor65} proved that
$d(z)^{n+1}K_\Om(z)$ has a positive finite limit at every
strictly pseudoconvex $C^2$ boundary point, where
$d(z)=\dist(z,\partial\Om)$.
At an $h$-extendible point of a smoothly bounded pseudoconvex domain
with Catlin multitype
$(1,m_2,\ldots,m_n)$, the nontangential growth exponent is
$2+\sum_{j=2}^n2/m_j$; this includes convex finite-type points
(Boas, Straube, and Yu \cite{BSY95}; see also
Chen and Fu \cite[Section~4, p.~77]{CF}).

For bounded pseudoconvex domains with Lipschitz boundary near the point,
Nikolov and Thomas \cite[Proposition~5]{NTkernel} give the lower bound
$K_\Om(z)\ge c\,d(z)^{-2}$ for $z$ near that point and some $c>0$.
They attribute the local $C^2$ case to Fu \cite[Theorem, p.~979]{Fu94}.
On a bounded pseudoconvex domain with smooth boundary,
Chen and Fu \cite[Theorem~4.2, p.~77]{CF} show that strict
pseudoconvexity holds at any point where, for some
$\beta>n+\tfrac12$, the nontangential lower limit of
$d(z)^\beta K_\Om(z)$ is positive; see their definition (4.2).
By contrast, \Cref{thm:normal-log-growth} obtains the exponent
$n+1$ from the full normal squeezing limit under only local
$C^2$ regularity, with no finite-type or convexity assumption.
Its conclusion is a logarithmic limit along the normal, not a
positive limit for $\delta^{n+1}K_\Om(p-\delta\nu_p)$.

\subsection{Proof strategy and organization}

A boundary curve germ would give a nonconstant boundary-valued limit
of the inverse charts $G_j=f_j^{-1}$. A holomorphic covector
proportional to the conormal by a positive real factor gives a
local retraction. A local kernel comparison, with a constant independent of the
transverse radius, then gives a contradiction by the minimum
principle for superharmonic functions. With
boundary curves excluded, $G_j\to p$
(\Cref{thm:no-curves,cor:inverse-collapse}).

To find the kernel growth exponent, we select normal scales using
kernel bounds and conjugate the rescaled normal translations by $f_j$.
A further subsequence yields a holomorphic semigroup. Its orbit
and Jacobian identify the normalized kernel limit at the selected
scales. If the lower limit of the logarithmic quotient in
\eqref{eq:normal-log-growth} were below $n+1$, scale selection would
give a lower bound at $s=2$ incompatible with that limit. Together
with the kernel upper bound, this gives the exponent $n+1$ along
the full normal. If a local defining function has a second-order
Taylor remainder of order $O(|h|^{2+\alpha})$ at $p$, Levi
degeneracy would give a kernel upper exponent strictly below
$n+1$ by \eqref{eq:weak-kernel-bound}. The contradiction proves
strict pseudoconvexity.

\Cref{sec:prelim,sec:boundary-images} give the preliminary estimates
and the exclusion of boundary curves. \Cref{sec:profiles,sec:classification}
identify the selected kernel limits, which give the main theorems
in \Cref{sec:conclusion}.
\Cref{sec:applications} treats intrinsic metrics, sampling, and
biholomorphic realizations.

\section{Preliminaries}\label{sec:prelim}

We fix notation and prove kernel-convergence and boundary-limit
lemmas for almost full ball charts. We then establish normal
Bergman-kernel bounds and select scales with useful lower bounds
for the normalized kernels.

\subsection{Notation and conventions}

A \emph{domain} is a nonempty connected open set. Unless explicitly
stated otherwise, $\Om\Subset\C^n$ is a bounded pseudoconvex
domain, $n\ge2$, and $p\in\partial\Om$ has a neighborhood in
which the boundary is $C^2$. No boundary regularity is assumed
elsewhere. For $0<\alpha\le1$, $C^{2,\alpha}$ means that the
second derivatives are locally H\"older continuous with exponent
$\alpha$. This regularity is used only to obtain the remainder bound
\eqref{eq:pointwise-taylor}, which improves the kernel upper exponent
at a Levi-degenerate point. The
normal squeezing condition \eqref{eq:normal-hypothesis} is imposed
only where stated.

For each integer $m\ge1$, write
$\BB^m=\{z\in\C^m:|z|<1\}$ and $\DD=\BB^1$.
We use $\C^0=\{0\}$, omit coordinate blocks of dimension zero,
and interpret empty sums as zero. In the relevant Euclidean space, let
$B(a,r)=\{z:|z-a|<r\}$ for $r>0$, and write $A\Subset U$
when $\overline A$ is compact and contained in $U$.
Norms are Euclidean norms or the induced operator norms.
For $\zeta\in\C$, $\Rea\zeta$ and $\Ima\zeta$ denote its
real and imaginary parts, and $\overline\zeta$ its complex conjugate.
Tangent vectors to open subsets of $\C^m$ are identified with
vectors in $\C^m$, with standard coordinate vectors $e_j$.
The Hermitian inner product
\[
 \ip{z}{w}=\sum_{j=1}^m z_j\overline{w_j}
 \qquad(z,w\in\C^m)
\]
is linear in the first variable. For a matrix $A$, $A^*$ denotes
its conjugate transpose; $\ker A$, $\ran A$, and $\rank A$ denote
its kernel, range, and rank, and $\tr A$ its trace when $A$ is square.
For Hermitian $A$, $A\ge0$ means positive semidefiniteness, and
$E^\perp$ is the Hermitian orthogonal complement of a complex
subspace $E$. Write $I$ for the identity linear map and $\id_X$
for the identity on $X$, suppressing the subscript when clear.

For open sets $U\subset\mathbb C^m$ and $V\subset\mathbb C^k$, let $\mathcal O(U,V)$ denote the set of holomorphic maps from $U$ to $V$. For domains of the same dimension, a
holomorphic embedding is an injective holomorphic map with a
holomorphic inverse on its open image.
In coordinates $z_b=x_b+iy_b$, the Wirtinger derivatives are
$\partial_{z_b}=\tfrac12(\partial_{x_b}-i\partial_{y_b})$ and
$\partial_{\overline z_b}=\tfrac12(\partial_{x_b}+i\partial_{y_b})$.
For a differentiable map, $d$ denotes its real differential.
If $m,k\ge1$, $U\subset\C^m$ is open, and
$f=(f_1,\ldots,f_k):U\to\C^k$ is holomorphic,
we identify its differential at $z\in U$ with the complex-linear map
\[
 df(z)=\left(\frac{\partial f_a}{\partial z_b}(z)\right)_{
 1\le a\le k,\,1\le b\le m}:\C^m\longrightarrow\C^k.
\]
We also write $df_z$ for $df(z)$; when $k=m$, $\det df(z)$ is
the complex Jacobian determinant. For a parameterized family of
holomorphic maps $\Phi(t,\cdot)$, write $\Phi_t(z)=\Phi(t,z)$
and use $d_z$ for the spatial differential with the parameter fixed,
so that $d\Phi_t(z)=d_z\Phi(t,z)$.

In the general definitions below, $D\subset\C^m$, $m\ge1$, is
an auxiliary domain, with no boundedness, pseudoconvexity, or
boundary regularity assumed unless stated; $\Aut(D)$ is its
biholomorphism group. The standing hypotheses on $\Om$ do not
apply to $D$. For Bergman kernels, $D$ may also be a nonempty
disconnected open set.
Constants $c,C>0$ may change from one occurrence to the next;
their uniformity is specified when needed.

\textbf{Boundary geometry.}
We write $\operatorname{PSH}(D)$ for the plurisubharmonic functions
on $D$: these are upper semicontinuous functions
$u:D\to[-\infty,\infty)$ whose restriction to each connected
component of the intersection of $D$ with a complex line is
subharmonic or identically $-\infty$. Pseudoconvexity means the
existence of a continuous plurisubharmonic exhaustion; for a
disconnected open set, it is understood componentwise.
Pluriharmonic functions are locally real parts of holomorphic functions.

For the boundary notions, suppose $D\subset\C^m$ has $C^2$
boundary near $q\in\partial D$. A local defining function is a
real-valued $C^2$ function $\rho$ on a neighborhood $U$ of $q$,
with $D\cap U=\{z\in U:\rho(z)<0\}$,
$\partial D\cap U=\{z\in U:\rho(z)=0\}$, and $d\rho\ne0$
on $\partial D\cap U$. Write $\rho_{z_k}=\partial_{z_k}\rho$,
and similarly for higher derivatives. The $(1,0)$ part of $d\rho$ is
$\partial\rho=\sum_{k=1}^m\rho_{z_k}\,dz_k$.
Subscripts such as $\partial\rho_q$ indicate evaluation at the
specified point. For $v\in\C^m$, viewed also as a real vector,
\[
 \partial\rho_q(v)=\sum_{k=1}^m\rho_{z_k}(q)v_k,
 \qquad d\rho_q(v)=2\Rea\bigl(\partial\rho_q(v)\bigr).
\]
The complex tangent space and Levi form are
\[
 T_q^{1,0}\partial D=\ker\partial\rho_q,\qquad
 \mathcal L_{\rho,q}(v,\overline w)
 =\sum_{k,l=1}^m\rho_{z_k\overline z_l}(q)v_k\overline{w_l}
 \quad(v,w\in\C^m).
\]
The boundary is \emph{strictly pseudoconvex}, or equivalently
\emph{strongly pseudoconvex}, at $q$ if this Levi form is positive
definite on $T_q^{1,0}\partial D$. For a pseudoconvex domain,
we call $q$ \emph{weakly pseudoconvex} when the restricted Levi
form has a nonzero null vector; its nullity is the \emph{Levi
corank}. In dimension one, positivity on the zero-dimensional
complex tangent space is vacuous. A conormal row at $q$ is a
nonzero complex-linear covector annihilating $T_q^{1,0}\partial D$;
its pairing with a column vector is written $\lambda v$, without
complex conjugation. The outward Euclidean unit normal is
$\nu_q=\nabla\rho(q)/|\nabla\rho(q)|$, so
$\partial\rho_q(\nu_q)=|\nabla\rho(q)|/2>0$.

A \emph{boundary curve germ} through $q\in\partial D$ is a
nonconstant holomorphic germ $\gamma:(\DD,0)\to(\C^m,q)$
with image in $\partial D$; its derivative may vanish at $0$.
For an open set $U\subset\C^m$, a holomorphic retraction is a
map $H\in\mathcal O(U,U)$ satisfying $H\circ H=H$; for a germ
fixing its base point, the identity is required only locally. The generic rank of a holomorphic map on a connected
source is its maximal rank, attained off the common zero set of its
maximal minors.

\textbf{Normal coordinates and convergence.}
Let $\nu_p$ be the outward Euclidean unit normal at $p$.
An affine unitary change of coordinates places $p$ at $0$ and
$\nu_p$ at $e_n$. We retain $p$ and
$\nu_p$ in formulas when they clarify the geometry. Set
\[
 \HH=\{s\in\C:\Rea s>0\},\qquad
 q_\delta(s)=p-\delta s\nu_p\quad(\delta>0,\ s\in\HH).
\]
Real $s>0$ traces the inward normal ray. The full domain $\HH$
parametrizes the inward half-plane in the complex normal line.
Write $s=x+iy$ with $x=\Rea s$ and $y=\Ima s$.
Local $C^2$ regularity gives a fixed interior tangent ball
$B(p-R\nu_p,R)\subset\Om$. For nonempty $Q\Subset\HH$,
$q_\delta(Q)$ lies in this ball when
$0<\delta<2R\inf_Q\Rea s/\sup_Q|s|^2$, since
$|R-\delta s|^2=R^2-2R\delta\Rea s+\delta^2|s|^2$.
Thus $q_\delta$ eventually sends every compact subset of $\HH$ into $\Om$.

Compact-open convergence means locally uniform convergence;
$C^\infty_{\mathrm{loc}}$ convergence includes all derivatives.
For holomorphic maps on varying domains, Montel's theorem is
applied on a fixed countable exhaustion of the limiting domain,
followed by a diagonal subsequence. We use the same convention
for maps on $r_j\BB^n$ when $r_j\to1$; the radii need not be
monotone.

\textbf{Invariant distances and metrics.}
For a domain $D\subset\C^m$, write
$d_D(z)=\dist(z,\partial D)$ for $z\in D$, with
$\dist(z,\varnothing)=+\infty$. The Kobayashi pseudodistance
$k_D$ is defined by chains of holomorphic discs equipped with
the Poincar\'e distance. We normalize it by
$k_{\DD}(0,t)=\operatorname{arctanh}t$ for $0\le t<1$.
The infinitesimal Kobayashi--Royden pseudometric, for
$z\in D$ and $v\in\C^m$, is
\[
 \kappa_D(z;v)=\inf\{a>0:f\in\mathcal O(\DD,D),\
                      f(0)=z,\ a f'(0)=v\}.
\]
In particular, these are a distance and a nondegenerate infinitesimal
metric on bounded domains and on $\HH$.
An approach $z\in D$ to a boundary point $q\in\partial D$ is
\emph{nontangential} if $|z-q|\le A d_D(z)$ for one fixed
$A\ge1$. A pseudohyperbolic ball of radius $r\in(0,1)$ in
$\BB^m$ is a $k_{\BB^m}$-ball of radius $\operatorname{arctanh}r$.

\textbf{Bergman kernels and metrics.}
For any nonempty open set $D\subset\C^m$, including a disconnected
one, let $A^2(D)$ be the space of square-integrable holomorphic
functions with respect to Euclidean volume $dV$; on a complex
line, $dA$ denotes Euclidean area. We denote its Bergman kernel, holomorphic in $z$ and
antiholomorphic in $w$, by $K_D(z,w)$ for $z,w\in D$, and
write $K_D(z)=K_D(z,z)$ for its diagonal.
On a disconnected open set, this is the kernel of the relevant
component when both arguments belong to that component, and
is zero when they belong to different components. In particular,
\[
 K_D(z)=\sup\bigl\{|f(z)|^2:f\in A^2(D),\
                          \textstyle\int_D|f|^2dV\le1\bigr\}.
\]
The diagonal kernel is positive on every bounded open set:
the nonzero constant function is square-integrable there.
It can vanish on an unbounded open set. The unbounded domains
for which we use Bergman kernels below are biholomorphic to the
ball, so their kernels are positive by the transformation law.

For nonempty open subsets of the same $\C^m$, we use the
monotonicity and biholomorphic transformation laws
\begin{align}
 D_1\subset D_2&\quad\Longrightarrow\quad
 K_{D_2}(z)\le K_{D_1}(z)\quad(z\in D_1),
 \label{eq:kernel-monotone}\\
 K_D(z)&=\abs{\det df(z)}^2K_{f(D)}(f(z))
 \quad\text{if }f:D\to f(D)\subset\C^m\text{ is biholomorphic}.
 \label{eq:kernel-transform}
\end{align}
They follow from the extremal definition and change of variables
in the $A^2$ norm. For the unit ball in $\C^m$,
\begin{equation}\label{eq:ball-kernel}
 K_{\BB^m}(z,w)=\frac{m!}{\pi^m}
 (1-\ip{z}{w})^{-(m+1)}\qquad(z,w\in\BB^m).
\end{equation}
On a bounded domain $D\subset\C^m$, the Bergman metric, for
$z\in D$ and $v\in\C^m$, is normalized by
\[
 b_D(z;v)^2=\sum_{a,b=1}^m
 \frac{\partial^2\log K_D(z)}{\partial z_a\partial\overline z_b}
 v_a\overline{v_b}.
\]
This metric is biholomorphically invariant: in the logarithm of
\eqref{eq:kernel-transform}, the Jacobian term is pluriharmonic
and hence makes no contribution to the complex Hessian.

\subsection{Almost full ball charts and their limits}

We call holomorphic embeddings $f_j:\Om\hookrightarrow\BB^n$
\emph{almost full ball charts} if there are $0<r_j<1$ with
\[
 r_j\BB^n\subset D_j:=f_j(\Om)\subset\BB^n,
 \qquad r_j\to1.
\]
The inverse $G_j=f_j^{-1}$ is defined on $D_j$. When taking
limits on $\BB^n$, we use its restriction to $r_j\BB^n$.
A chart is centered at $a_j\in\Om$ when $f_j(a_j)=0$.
For every sequence $\delta_j\downarrow0$, the points $q_{\delta_j}(1)$
lie in $\Om$ for all sufficiently large $j$. By the definition of
squeezing, \eqref{eq:normal-hypothesis} gives almost full ball charts
centered at these points. The next two lemmas do not require that limit.

\begin{lem}[Kernel convergence in almost full balls]\label{lem:ball-kernel}
Suppose $r_j\BB^n\subset D_j\subset\BB^n$ and $r_j\to1$,
where each $D_j$ is open. Then
$K_{D_j}(z,w)\to K_{\BB^n}(z,w)$ on compact subsets of
$\BB^n\times\BB^n$, with all derivatives in $z$ and $\overline w$.
\end{lem}

\begin{proof}
First control the diagonal. Monotonicity gives
$K_{\BB^n}(z)\le K_{D_j}(z)\le K_{r_j\BB^n}(z)$.
By \eqref{eq:kernel-transform} and \eqref{eq:ball-kernel},
$K_{r\BB^n}(z)=r^{-2n}K_{\BB^n}(z/r)$, so both bounds tend
locally uniformly to $K_{\BB^n}(z)$. Off the diagonal,
\[
 |K_{D_j}(z,w)|^2\le K_{D_j}(z)K_{D_j}(w).
\]
Thus the kernels form a locally bounded holomorphic family in
$(z,\overline w)$ on compact sets eventually in their domains.
Every subsequential limit agrees with $K_{\BB^n}$ on the diagonal,
which determines all holomorphic--antiholomorphic Taylor coefficients.
This identifies the limit, and Cauchy's estimates give convergence
with all derivatives.
\end{proof}

\begin{lem}[Boundary limits of inverse ball charts]\label{lem:boundary-limit}
Let $f_j:\Om\hookrightarrow\BB^n$ be almost full ball charts
on a bounded domain $\Om\subset\C^n$, with inverses $G_j$.
If $G_j\to G$ locally uniformly on $\BB^n$ and
$G(0)\in\partial\Om$, then $G(\BB^n)\subset\partial\Om$.
If no nonconstant holomorphic curve germ through $G(0)$ lies
in $\partial\Om$, then $G$ is constant.
\end{lem}

\begin{proof}
For fixed $z\in\BB^n$, we have $G_j(z)\in\Om$ for all large $j$,
so $G(z)\in\overline\Om$. We exclude interior values by showing
that one would force an inverse chart to hit the boundary.

Suppose $q=G(a)\in\Om$. Then $a\ne0$. Since $|f_j|\le1$,
Cauchy's estimates on a ball about $q$ compactly contained in
$\Om$ bound $\|df_j(G_j(a))\|$ by $M$ once $G_j(a)$ lies in its
concentric half-ball.
Differentiating $f_j\circ G_j=\id$ at $a$ gives
$|v|\le M|dG_j(a)v|$. Local uniform convergence of the derivatives
makes $dG(a)$ invertible. On a disc in the complex line through
$0,a$ containing both points and compactly contained in $\BB^n$,
the determinants $\det dG_j$ eventually never vanish. Their
limit is nonzero at $a$, so Hurwitz's theorem gives $\det dG(0)\ne0$.

Set $p_0=G(0)$ and $A=dG(0)$. Choose $0<b<1$ with
$\|I-A^{-1}dG\|<1/4$ on $\overline{b\BB^n}$. For large $j$ this ball
lies in $r_j\BB^n$. Since $dG_j\to dG$ there and $G_j(0)\to p_0$, the map
\[
 T_j(z)=z-A^{-1}(G_j(z)-p_0)
\]
has $\|dT_j\|\le1/2$ on this ball and $|T_j(0)|<b/2$. It is
therefore a contraction of the closed ball into itself. Its
fixed point $z_j$ satisfies $G_j(z_j)=p_0$, impossible since
$G_j(z_j)\in\Om$ and $p_0\in\partial\Om$. Thus $G(\BB^n)\subset\partial\Om$.

If $G$ is nonconstant, choose $w$ with $G(w)\ne G(0)$.
The map $\zeta\mapsto G(\zeta w)$ is holomorphic on a disc
containing $0$ and $1$. By the identity theorem its germ at
$0$ is a nonconstant boundary curve through $G(0)$.
\end{proof}

The inverse-chart argument needs neither pseudoconvexity nor
boundary regularity. An alternative argument, on smoothly
bounded pseudoconvex domains, uses the negative strictly
plurisubharmonic exhaustion constructed by Diederich and
Forn\ae ss \cite[Theorem~I, p.~129]{DF}.

\subsection{Normal Bergman-kernel bounds}

Under the standing boundedness and pseudoconvexity assumptions,
the rough two-sided estimate requires only local $C^2$ regularity.
At a Levi-degenerate point, an order-$2+\alpha$ bound on the
Taylor remainder improves the kernel upper exponent to a number
strictly below $n+1$. At the same Levi-degenerate point, $C^2$
regularity alone gives $K_\Om(q_\delta(1))=o(\delta^{-(n+1)})$.

\begin{lem}\label{lem:kernel-bounds}
Let $\Om\Subset\C^n$ be a pseudoconvex domain, and suppose that its
boundary is $C^2$ near $p\in\partial\Om$.
For all sufficiently small $\delta>0$,
\begin{equation}\label{eq:rough-kernel-bounds}
 c\delta^{-2}\le K_\Om(q_\delta(1))
 \le C\delta^{-(n+1)}.
\end{equation}
Suppose that the tangential Levi form at $p$ has
$\ell\in\{1,\ldots,n-1\}$ zero eigenvalues.
If the boundary is $C^{2,\alpha}$ near $p$, $0<\alpha\le1$, then
\begin{equation}\label{eq:weak-kernel-bound}
 K_\Om(q_\delta(1))
 \le C\delta^{-a_{\alpha,\ell}},
 \qquad
 a_{\alpha,\ell}=n+1-\frac{\ell\alpha}{2+\alpha}<n+1.
\end{equation}
At the same Levi-degenerate point, with only $C^2$ regularity,
one still has
\begin{equation}\label{eq:weak-kernel-little-o}
 \delta^{n+1}K_\Om(q_\delta(1))\longrightarrow0.
\end{equation}
The power bound \eqref{eq:weak-kernel-bound} also holds
if, instead of local H\"older regularity, some $C^2$
defining function $\rho$ near $p$ satisfies
\begin{equation}\label{eq:pointwise-taylor}
 |\rho(p+h)-T_p^2\rho(h)|\le C|h|^{2+\alpha},
 \qquad 0<\alpha\le1,
\end{equation}
for all sufficiently small real increments $h$, where
$T_p^2\rho$ is its real Taylor polynomial of degree two at $p$.
The constants may depend on $\Om$ and $p$.
\end{lem}

\begin{proof}
We first obtain the rough bounds at the exact normal depth, then
use the Levi-null directions to enlarge one comparison polydisc.
The fixed tangent ball above contains $B(q_\delta(1),\delta)$
for $0<\delta<R$. Since $p\in\partial\Om$, the depth is $\delta$.

\proofstep{1}{The two-sided bound}
For bounded pseudoconvex domains with boundary Lipschitz near $p$,
Nikolov and Thomas \cite[Proposition~5]{NTkernel} prove
$\liminf_{\Om\ni z\to p}d_\Om(z)^2K_\Om(z)>0$.
The local $C^2$ graph meets this hypothesis without any regularity
away from $p$. Evaluating at depth $\delta$ gives the lower bound in
\eqref{eq:rough-kernel-bounds}; see also the local $C^2$ result of
Fu \cite[Theorem, p.~979]{Fu94}. For the upper bound, scaling \eqref{eq:ball-kernel} gives the tangent-ball kernel at $q_\delta(1)$ as $n!R^2\pi^{-n}(2R\delta-\delta^2)^{-(n+1)}$.
Kernel monotonicity gives the upper bound in
\eqref{eq:rough-kernel-bounds}.

\proofstep{2}{Larger radii in the Levi-null directions}
In affine unitary normal coordinates $(z,w)\in\C^{n-1}\times\C$,
diagonalize the tangential Levi form, semipositive by pseudoconvexity,
and remove the pure holomorphic tangential quadratic terms of $\rho$
via $(z,w)\mapsto(Uz,w+Q(z))$ with $U$ unitary and $Q$ homogeneous quadratic. The map fixes the normal line pointwise
and its Jacobian has modulus one. After positive normalization,
the defining function satisfies
\[
 \rho(z,w)=\Rea w+
 \sum_{k=\ell+1}^{n-1}\lambda_k|z_k|^2
 +O(|z||w|+|w|^2)+E(z,w),\qquad \lambda_k>0.
\]
Put $t=|z|+|w|$.
The shear adds $O(|z|^2|w|)$ to the mixed error and $O(t^3)$
to the Taylor remainder, whose order is preserved by the locally
bi-Lipschitz change. Thus $E=o(t^2)$ in the $C^2$ case and
$E=O(t^{2+\alpha})$ when the original defining function has an
order-$2+\alpha$ Taylor remainder at $p$, since $0<\alpha\le1$.

Fix $0<c\le1/4$ so small that
$c^2\sum_{k>\ell}\lambda_k\le1/4$.
For a parameter $r_\delta\to0$, take the polydisc with center
$(0,-\delta)$ and radii
\[
 r_\delta\ (1\le k\le\ell),\qquad
 c\sqrt\delta\ (\ell<k<n),\qquad c\delta\ (w).
\]
On this polydisc the sum of the normal linear term and the
Levi quadratic is at most $-\delta/2$. The mixed and normal
error is $O(\delta(r_\delta+\sqrt\delta)+\delta^2)=o(\delta)$
uniformly there. The polydisc shrinks into the defining neighborhood,
where $|E|\le\delta/4$ then yields $\rho<0$.
The reciprocal volume is the center kernel. Monotonicity and the
unit-modulus Jacobian give
\begin{equation}\label{eq:weak-polydisc-kernel}
 K_\Om(q_\delta(1))
 \le \pi^{-n}c^{-2(n-\ell)}\,
          \delta^{-(n+1-\ell)}r_\delta^{-2\ell}.
\end{equation}

When $E=O(t^{2+\alpha})$, set $r_\delta=c\delta^{1/(2+\alpha)}$. For each fixed $c$, take
$\delta$ small enough that $(1+c)\delta\le 2r_\delta$.
Then $t=O(r_\delta)$ on the polydisc with a constant
independent of $c$, so $|E|=O(c^{2+\alpha}\delta)$. Shrink $c$ once if needed to
make this at most $\delta/4$. Substitution in
\eqref{eq:weak-polydisc-kernel} yields \eqref{eq:weak-kernel-bound}.

Under $C^2$ regularity, keep $c$ fixed and take
$r_\delta=A\sqrt\delta$ for any fixed $A>0$.
Now $t=O_A(\sqrt\delta)$, so $E=o_A(\delta)$, and the
polydisc is included for sufficiently small $\delta$,
with a threshold allowed to depend on $A$. Thus the same estimate gives
\[
 \limsup_{\delta\downarrow0}
     \delta^{n+1}K_\Om(q_\delta(1))
 \le \pi^{-n}c^{-2(n-\ell)} A^{-2\ell}.
\]
Letting $A\to\infty$ proves \eqref{eq:weak-kernel-little-o}.
\end{proof}

\subsection{Scale selection and normalized Jacobians}

We select scales that give lower bounds for the normalized kernel at
small $x$ and on intervals extending farther along the positive real axis.

\begin{lem}[Two lower bounds with different slopes]\label{lem:envelope}
Let $h:[t_0,\infty)\to\R$ be continuous and satisfy
$h(t)\ge bt-C$ for some $b>0$. Suppose $B>b$ and
\begin{equation}\label{eq:envelope-subcritical}
 \liminf_{t\to\infty}\bigl(h(t)-Bt\bigr)=-\infty.
\end{equation}
Set $\mu=b/2$. There is a sequence $T_j\to\infty$ such that
\begin{equation}\label{eq:asymmetric-contact}
 h(t)-h(T_j)\ge
 \begin{cases}
  \mu(t-T_j),&t\ge T_j,\\
  B(t-T_j),&t_0\le t\le T_j.
 \end{cases}
\end{equation}
\end{lem}

\begin{proof}
To obtain both slopes from one minimizing point, choose $u_j\to\infty$
with $h(u_j)-Bu_j\to-\infty$ by \eqref{eq:envelope-subcritical}, and set
\[
 \ell_j(t)=\min\{Bt,\,\mu t+(B-\mu)u_j\}.
\]
For each fixed $j$, the right-hand slope $\mu<b$ makes $h-\ell_j$
tend to $+\infty$, so it attains a minimum at $T_j\in[t_0,\infty)$.
The minimum values tend to $-\infty$ by evaluation at $u_j$.
On each fixed compact interval $h-\ell_j=h(t)-Bt$ eventually,
so $T_j\to\infty$.

At the minimum, $h(t)-h(T_j)\ge\ell_j(t)-\ell_j(T_j)$.
All secant slopes of $\ell_j$ lie between $\mu$ and $B$.
Use the lower slope bound for $t\ge T_j$ and the upper slope
bound for $t\le T_j$, where $t-T_j\le0$, to obtain
\eqref{eq:asymmetric-contact}.
Take a subsequence with $T_j$ strictly increasing.
\end{proof}

\begin{remark}
The familiar upper bound $h(t)\le at+C$ with $a<B$ implies
\eqref{eq:envelope-subcritical}. So does the weaker condition
$\liminf_{t\to\infty}h(t)/t<B$. In the proof, $b/2$ may be
replaced by any fixed number in $(0,b)$; the stated choice is
sufficient here.
\end{remark}

To apply the lemma to the normal kernel, choose $t_0$ large
enough that $q_{e^{-t}}(1)\in\Om$ for $t\ge t_0$, and set
\[
 h(t)=\log K_\Om(q_{e^{-t}}(1)).
\]
The rough bounds in \Cref{lem:kernel-bounds} give
$2t-C\le h(t)\le(n+1)t+C$. Thus we take $b=2$ and
$\mu=1$; each application specifies a $B>2$ satisfying
\eqref{eq:envelope-subcritical}. For the resulting scales
$\delta_j=e^{-T_j}$ and $0<x\le e^{T_j-t_0}$, put
$R_j(x)=K_\Om(q_{\delta_j}(x))/K_\Om(q_{\delta_j}(1))$.
Then $\log R_j(x)=h(T_j-\log x)-h(T_j)$.
In \eqref{eq:asymmetric-contact}, $x\le1$ uses $t\ge T_j$,
while $x\ge1$ uses $t\le T_j$. Exponentiating gives
\begin{equation}\label{eq:asymmetric-profile-bound}
 R_j(x)\ge
 \begin{cases}
 x^{-\mu},&0<x\le1,\\
 x^{-B},&1\le x\le e^{T_j-t_0}.
 \end{cases}
\end{equation}
Since $T_j\to\infty$, every compact interval $I\subset(0,\infty)$
eventually lies below $e^{T_j-t_0}$. On $I$, the two bounds give
$R_j(x)\ge\max\{1,\sup I\}^{-B}$ for all such $j$.
The bound $k(x)\ge x^{-\mu}$ forces a normalized kernel limit $k$
to diverge as $x\downarrow0$, while $k(x)\ge x^{-B}$ for $x\ge1$
prevents exponential decay. The first bound is used to exclude
boundary curves; both enter the semigroup argument. We select
scales from $h$, then use \eqref{eq:normal-hypothesis} to choose
the charts.

The bounds in \eqref{eq:asymmetric-profile-bound} concern only
the positive real axis. The next lemma gives compactness of the
normalized Jacobians once kernel ratios are bounded away from
zero on compact subsets of $\HH$. The kernel comparisons below
supply this bound.

\begin{lem}\label{lem:jacobian-compactness}
Let $\delta_j\downarrow0$ and let $f_j$ be almost full ball
charts centered at $q_{\delta_j}(1)$. After passage to a subsequence, write
$F_j(s)=f_j(q_{\delta_j}(s))\to F(s)$ locally uniformly on
$\HH$. Suppose that for every $Q\Subset\HH$ there is $c_Q>0$
such that
\begin{equation}\label{eq:profile-positive}
 R_j(s):=\frac{K_\Om(q_{\delta_j}(s))}
                    {K_\Om(q_{\delta_j}(1))}\ge c_Q
 \qquad(s\in Q)
\end{equation}
for all sufficiently large $j$. Then a further subsequence
satisfies
\begin{equation}\label{eq:jacobian-ratio}
 J_j(s):=\frac{\det df_j(q_{\delta_j}(s))}
                    {\det df_j(q_{\delta_j}(1))}
 \longrightarrow J(s),
\end{equation}
locally uniformly on $\HH$, where $J$ is holomorphic and
nowhere zero and $J(1)=1$. The normalized kernels converge in
$C^\infty_{\mathrm{loc}}(\HH)$:
\begin{equation}\label{eq:profile-formula}
 R_j(s)\longrightarrow
 k(s):=\frac{|J(s)|^2}{(1-|F(s)|^2)^{n+1}}.
\end{equation}
\end{lem}

\begin{proof}
Put $D_j=f_j(\Om)$. The kernels at the chart centers are positive
by boundedness, and each $J_j$ is nonvanishing because $f_j$ is
biholomorphic onto $D_j$. Since $F(1)=0$, the maximum principle puts
$F$ in the open ball. The transformation law \eqref{eq:kernel-transform} gives
\begin{equation}\label{eq:kernel-jacobian}
 R_j(s)=|J_j(s)|^2A_j(s),\qquad
 A_j(s)=\frac{K_{D_j}(F_j(s))}{K_{D_j}(0)}.
\end{equation}
On each compact parameter set the images under $F_j$ lie eventually in
a fixed compact subball. By \Cref{lem:ball-kernel} and Cauchy's estimates
for $F_j$, the chain rule gives $A_j\to(1-|F|^2)^{-(n+1)}$
in $C^\infty_{\mathrm{loc}}(\HH)$.
Since $|1/J_j|^2=A_j/R_j$, \eqref{eq:profile-positive}
bounds the holomorphic reciprocals locally. Montel gives a
subsequential limit $g$, and $g(1)=1$ rules out zeros by Hurwitz.
Thus $J_j\to1/g$ locally uniformly with all derivatives.
Multiplying $|J_j|^2$ by $A_j$ gives \eqref{eq:profile-formula}
in $C^\infty_{\mathrm{loc}}(\HH)$.
\end{proof}

\section{Boundary geometry and limits of inverse charts}\label{sec:boundary-images}

Fix a sequence of almost full ball charts. Throughout this section,
only the boundary near $p$ is assumed to be $C^2$. Let $\rho$ be a defining
function on this neighborhood. For a transverse vector $v$
with
\begin{equation}\label{eq:transverse-vector}
 \partial\rho_p(v)>0\quad\text{real},
\end{equation}
suppose
\begin{equation}\label{eq:transverse-charts}
 f_j(p-\delta_jv)=0,\qquad
 r_j\BB^n\subset f_j(\Om)\subset\BB^n,\qquad
 \delta_j\downarrow0,\quad r_j\to1.
\end{equation}
By boundedness and \Cref{lem:boundary-limit}, after passing to a subsequence we have
\begin{equation}\label{eq:boundary-G}
 G_j:=f_j^{-1}\longrightarrow G:\BB^n\to\partial\Om,
 \qquad G(0)=p.
\end{equation}
A holomorphic conormal along $G$ near $0$ gives a local
retraction and hence constant rank. The localization estimate
also uses boundedness of $G$ on the whole ball, but it needs
boundary regularity only near $p$.

\subsection{Boundary curves and positive holomorphic conormals}

We first show that a boundary curve forces $dG(0)$ to have
positive rank. We then construct a positive holomorphic conormal
along the local boundary image, without assuming constant rank.

\begin{lem}\label{lem:curve-rank}
In the setting of \eqref{eq:transverse-charts}--\eqref{eq:boundary-G}, suppose
$\gamma(\zeta)=p+A\zeta^m+O(\zeta^{m+1})$ is a boundary curve
germ, with integer $m\ge1$ and $A\in\C^n\setminus\{0\}$. Then $A\in\ran dG(0)$. In particular,
$\rank dG(0)>0$ whenever a boundary curve germ through $p$ exists.
\end{lem}

\begin{proof}
We recover $A$ from the first nonzero coefficient of the translated
curve. On a fixed small closed source disc about $0$ in the parameter
domain of $\gamma$, $\Rea\partial\rho_{\gamma(\zeta)}(v)>c>0$.
The derivative of $\rho$ in the direction $-v$ is less than
$-2c$. Hence $\gamma(\zeta)-\delta_jv$ lies in $\Om$ for large
$j$, uniformly on that disc. The holomorphic maps
\[
 H_j(\zeta)=f_j(\gamma(\zeta)-\delta_jv)
\]
take values in $\BB^n$ and vanish at $0$. Since $\gamma-p$ has no
terms of degree below $m$, the coefficient of $\zeta^m$ in $H_j$
is $df_j(p-\delta_jv)A$. Cauchy's estimate on the fixed disc
bounds these coefficients, so a subsequence converges to $b$.
The inverse identity $dG_j(0)\,df_j(p-\delta_jv)A=A$ gives
$dG(0)b=A$ in the limit.
\end{proof}

We need a holomorphic inward field along $G$, with a first
derivative that annihilates $\ker dG$. To obtain it, first use
Levi semipositivity to show that dividing by a nonzero component gives a holomorphic row.
The remaining issue is to make the multiplier positive real.
We prove this by moving transverse analytic discs and using
pseudoconvexity to control their centers. The argument does not
need $G$ to have constant rank. For related work, see Bedford and
Forn\ae ss \cite{BF81} on complex manifolds in pseudoconvex
boundaries and Baracco \cite{Bar03} on analytic discs in
conormal bundles.

\begin{lem}[A local positive holomorphic conormal]\label{lem:conormal}
Let $\Om\subset\C^n$ be a pseudoconvex domain and let
$G:(\C^m,0)\to(\partial\Om,p)$, $m\ge1$, be holomorphic.
Let $\rho$ be a $C^2$ defining function near $p$, oriented so
that $\Om$ is locally $\{\rho<0\}$, and choose $v\in\C^n$
with $\partial\rho_p(v)>0$ real.
On a sufficiently small source ball $U$, there are a positive
$C^1$ function $\beta$ and a nowhere-zero holomorphic row $\lambda$
whose holomorphic field $V=v/(\lambda v)$ satisfies
\begin{align}
 \lambda&=\beta(\partial\rho\circ G),\qquad \lambda(0)v=1,
 \label{eq:positive-conormal}\\
 \lambda V&=1,\qquad V(0)=v,\qquad dV(0)|_{\ker dG(0)}=0,
 \label{eq:normal-field-jet}\\
 d\lambda(z)|_{\ker dG(z)}&=0\qquad(z\in U).
 \label{eq:conormal-kernel-jet}
\end{align}
For every $K\Subset U$ and $T\Subset\HH$, and all sufficiently small $\delta>0$,
\begin{equation}\label{eq:field-inward}
 G(z)-\delta\tau V(z)\in\Om
 \qquad(z\in K,\ \tau\in T).
\end{equation}
\end{lem}

\begin{proof}
We first find a holomorphic row in the conormal direction, then
use a disc test to make its factor positive.

\proofstep{1}{A holomorphic conormal}\par\nopagebreak[4]
\noindent Write $\vartheta=\partial\rho\circ G$. Choose a coordinate
vector $e_k$ so that $c=\vartheta(e_k)$ is nowhere zero after
shrinking $U$, and take a $C^1$ branch $\psi=\arg c$.
Differentiating $\rho\circ G=0$ puts $dG$ in the complex tangent space.
Since $G$ is holomorphic, $\partial\overline\partial(\rho\circ G)=0$
makes its image Levi-null. For a source direction $v_0$ the derivative
$\overline\partial_{v_0}\vartheta$ is the row
$\xi\mapsto\mathcal L_{\rho,G(z)}(\xi,\overline{dG(z)v_0})$. By semipositivity,
the row annihilates the complex tangent space and equals
$(\overline\partial_{v_0}c)\vartheta/c$ by evaluation at $e_k$.
Thus $\lambda_0=\vartheta/c$ is holomorphic.

\proofstep{2}{Pluriharmonicity of the phase}\par\nopagebreak[4]
\noindent Center a small closed
complex-line disc $\Delta\Subset U$ at $z_0$ on which
$\chi=\pi-\psi$ oscillates by less than $\pi/8$. For
\[
 B_\varepsilon(z,w)=G(z)+\varepsilon e^w e_k
\]
Taylor's formula gives
$\rho(B_\varepsilon(z,w))/\varepsilon\to2\Rea(c(z)e^w)$
uniformly on compact sets in $\Delta\times\C$.
Since $\psi+\Ima w=\pi+(\Ima w-\chi)$, this limit is negative
in the central strip $|\Ima w-\chi(z)|<\pi/2$.

Let $H_\chi$ be the harmonic extension of $\chi|_{\partial\Delta}$.
For $0<\eta<\pi/8$, approximate the boundary values by a trigonometric
polynomial and extend harmonically to obtain a polynomial $H$
with $|H-H_\chi|<\eta/4$ on $\Delta$ by the maximum principle.
Choose a holomorphic polynomial $W$ with $\Ima W=H$ and move the
same graph in both directions, $\sigma\in\{-1,1\}$:
\[
 w_t^\sigma(z)=W(z)+i\sigma t(\pi/2-\eta),\qquad 0\le t\le1.
\]
The starting graph lies in the strip on $\Delta$ because
$|H-\chi|<\eta/4+\pi/8$ by the maximum principle and oscillation bound.
All boundary graphs also lie there: $|H-\chi|<\eta/4$ on the circle
and $|t(\pi/2-\eta)|\le\pi/2-\eta$. These strict margins and uniform
Taylor convergence put the starting disc and all boundary circles
in $\Om$ under $B_\varepsilon$ for small $\varepsilon>0$.

For each such $\varepsilon$, a continuous plurisubharmonic
exhaustion is bounded above by $M$ on the compact boundary-circle
images. The maximum principle puts every included closed disc
in the compact closed sublevel at height $M$. For either sign,
the inclusion parameters contain $0$, are open by compactness
of the disc, and are closed by the sublevel estimate. Hence $t=1$
also gives inclusion.

As $\varepsilon\downarrow0$, the two endpoints at the center
give $\Rea(c(z_0)e^{w_1^\sigma(z_0)})\le0$. The angles
$H(z_0)-\chi(z_0)+\sigma(\pi/2-\eta)$ belong to $(-\pi,\pi)$;
the weak sign puts both in $[-\pi/2,\pi/2]$, giving
$|H(z_0)-\chi(z_0)|\le\eta$. Letting $\eta\downarrow0$
yields $H_\chi(z_0)=\chi(z_0)$. The mean-value identity on
every such disc makes $\psi$ harmonic on complex lines,
hence pluriharmonic.

\proofstep{3}{Positive normalization and an inward field}\par\nopagebreak[4]
\noindent Choose holomorphic $q$ with $\Ima q=\psi$, shrinking $U$ if needed.
Then
\[
 \lambda=e^q\lambda_0=\beta\vartheta,\qquad
 \beta=\frac{e^{\Rea q}}{|c|}>0
\]
is a holomorphic nowhere-zero row with a positive $C^1$ factor.
Multiplying by a positive constant gives the normalization $\lambda(0)v=1$.

If a real tangent vector $\xi$ represents a vector in
$\ker dG(z)$, then $d\vartheta(\xi)=d\vartheta(i\xi)=0$. Hence
\[
 d\lambda(\xi)=d(\log\beta)(\xi)\lambda,
 \qquad d\lambda(i\xi)=d(\log\beta)(i\xi)\lambda.
\]
Holomorphy gives $d\lambda(i\xi)=i\,d\lambda(\xi)$.
The two displayed scalar coefficients are real, so both
vanish. This proves \eqref{eq:conormal-kernel-jet}.
Shrink $U$ so that $\lambda v\ne0$ and set $V=v/(\lambda v)$.
Then $V(0)=v$ and $dV=-v\,d(\lambda v)/(\lambda v)^2$;
\eqref{eq:conormal-kernel-jet} gives $dV|_{\ker dG}=0$.
Finally, uniform $C^2$ Taylor expansion gives
\[
 \rho(G(z)-\delta\tau V(z))
 =-2\delta\beta(z)^{-1}\Rea\tau+O_{K,T}(\delta^2)<0
\]
on the stated compact sets, proving \eqref{eq:field-inward}.
\end{proof}

\subsection{Local retractions and holomorphic tubes}

We do not initially know whether the image germ of $G$ is a
complex submanifold. We first construct a holomorphic retraction of the source
through which $G$ factors. This yields constant rank near the
origin. When this rank is positive, the image is a local complex
submanifold to which the holomorphic tube construction applies.

\begin{prop}[Local retraction and the generic rank]\label{prop:local-rank}
Let $\Om\subset\C^n$ be a bounded pseudoconvex domain with $C^2$
boundary near $p$. Let $f_j$ be almost full ball charts centered
at $p-\delta_jv$, where $\delta_j\downarrow0$ and
$\partial\rho_p(v)>0$ is real for a local defining function $\rho$.
Suppose $G_j=f_j^{-1}\to G$ locally uniformly on
$\BB^n$. Put $L=dG(0)$ and $r=\rank L$.
Then $G$ has generic rank $r$ on $\BB^n$ and constant rank
$r$ near $0$. More precisely, there is a holomorphic
retraction germ $H:(\C^n,0)\to(\C^n,0)$ of rank $r$
such that $G\circ H=G$ as germs at $0$.
\end{prop}

\begin{proof}
\proofstep{1}{A map preserving $G$} We construct a retraction onto an $r$-dimensional germ.
\Cref{lem:boundary-limit} places $G(\BB^n)$ in $\partial\Om$,
so \Cref{lem:conormal} supplies on a small source ball $U$
$\lambda=\beta(\partial\rho\circ G)$, $\beta>0$, and
$V=v/(\lambda v)$, with $V(0)=v$ and
$dV(0)|_{\ker L}=0$. By the inward-field assertion, the maps
\[
 H_j(z)=f_j(G(z)-\delta_jV(z))
\]
are eventually defined on every compact subset of $U$, take values
in $\BB^n$, and fix $0$. A subsequential limit $H$ takes values in
the open ball by the maximum principle at $H(0)=0$. Passing
to the limit gives
\[
 G\circ H=G\qquad\text{on }U.
\]
Indeed, for each compact subset of $U$, its images under $H_j$
eventually lie in a fixed compact subset of $\BB^n$,
where $G_j\circ H_j=G-\delta_jV$.

Set $P=dH(0)$. Since $H(0)=0$ and $G\circ H=G$, we have
$LP=L$. For each sufficiently large $j$, the chain rule for
$H_j=f_j\circ(G-\delta_jV)$ and $dV(0)|_{\ker L}=0$ show that
$dH_j(0)$ vanishes on $\ker L$. Passing to the limit gives
$\ker L\subset\ker P$, and $LP=L$ gives the reverse inclusion.
Hence $\rank P=r$.

\proofstep{2}{Idempotence and rank of $H$} Take a connected neighborhood
$W\ni0$ with $H(W)\subset U$. Since $G\circ H=G$,
the holomorphic quotient $(\lambda v)\circ H/(\lambda v)$
equals $(\beta\circ H)/\beta>0$ on $W$.
It is therefore constant, equal to one at $0$, and $V\circ H=V$.
For every compact subset of $W$ and all large $j$ this gives
the exact identity
\[
 H_j(H(z))
 =f_j(G(H(z))-\delta_jV(H(z)))
 =f_j(G(z)-\delta_jV(z))=H_j(z).
\]
Passing to the limit gives $H\circ H=H$ as germs at $0$,
and differentiation gives $P^2=P$.

To read off the rank, use the splitting $\C^n=\ran P\oplus\ker P$
and the coordinate change
\[
 z\longmapsto
 \bigl(PH(z),(I-P)(z-H(z))\bigr).
\]
Its derivative at zero is the identity in this splitting;
in these coordinates $H$ is $(x,y)\mapsto(x,0)$.
Since $G\circ H=G$, the rank of $G$ is at most $r$ near
$0$. All $(r+1)$-minors of $dG$ therefore vanish identically
on the connected ball. Since $\rank dG(0)=r$, lower
semicontinuity of rank gives constant rank $r$ near $0$,
and the generic rank is exactly $r$.
\end{proof}

Call the nonsingular image germ the \emph{local boundary leaf};
this term does not assert a boundary foliation. The next coordinates
make it the zero section while preserving the original complex
normal line.

\begin{lem}[A local holomorphic tube]\label{lem:tube}
Let $\Om\subset\C^n$ be a pseudoconvex domain, with $C^2$ boundary near
$p$, and let $M$ be a nonsingular complex submanifold germ of
dimension $r$, with $1\le r\le n-1$, contained in that boundary germ.
There are local holomorphic coordinates
$(\eta,y,w)\in\C^r\times\C^{n-r-1}\times\C$ in which
\begin{equation}\label{eq:tube-map}
 \Psi(\eta,y,w)=g(\eta)+\sum_{l=1}^{n-r-1}Y_l(\eta)y_l
                         +wV_0(\eta),
\end{equation}
where $g$ parametrizes $M$, $g(0)=p$, and
\begin{equation}\label{eq:tube-normal}
 \Psi(0,0,w)=p+w\nu_p.
\end{equation}
Writing $t=\Ima w$, the domain is locally defined by
\begin{equation}\label{eq:tube-graph}
 \Rea w+\phi(\eta,y,\Ima w)<0,
\end{equation}
where $\phi$ is $C^2$ and
\[
 \phi(\eta,0,0)=0,\qquad d_{y,t}\phi(\eta,0,0)=0.
\]
Here $d_{y,t}$ is the real differential in the fiber variables $(y,t)$.
For every $\varepsilon>0$ and any fixed smaller compact base, the bound
$|\partial_t\phi|<\varepsilon$ holds after shrinking only
the fiber variables $(y,w)$.
\end{lem}

\begin{proof}
The tube requires a holomorphic frame transverse to the leaf.
Apply \Cref{lem:conormal} to $g$ with $v=\nu_p$, so
$\lambda\,dg=0$ and $V_0=\nu_p/(\lambda\nu_p)$ satisfies
$\lambda V_0=1$ and $V_0(0)=\nu_p$.
Choose constant vectors $u_l$ completing the image of $dg(0)$
in $\ker\lambda(0)$, and set $Y_l=u_l-\lambda(u_l)V_0$.
Then $\lambda Y_l=0$ and $Y_l(0)=u_l$, so
$[dg,Y_1,\ldots,Y_{n-r-1},V_0]$ has full rank after shrinking.
The inverse function theorem gives injectivity of \eqref{eq:tube-map}
on a small product neighborhood. On the zero section,
$d(\rho\circ\Psi)$ annihilates the $\eta$ and $y$ directions,
and its $w$ part is a positive multiple of $\Rea dw$.
The real implicit function theorem gives \eqref{eq:tube-graph}
and the first-derivative vanishing; continuity of $\partial_t\phi$
gives the last assertion uniformly on a compact base.
\end{proof}

\subsection{Fixed-constant localization in thin tubes}

Under the limiting inverse chart $G$, pseudohyperbolic balls in $\BB^n$
centered at the limiting chart values $F(s)$ with a fixed small radius map into
a single compact piece of the local leaf. This gives a kernel
comparison independent of tube thickness, before $G_j\to p$ is known.

\begin{lem}[Uniform localization in a local boundary leaf]
\label{lem:local-leaf-localization}
Suppose $\Om\Subset\C^n$ is a bounded domain,
$p\in\partial\Om$, $v\in\C^n\setminus\{0\}$ is fixed, and
$\delta_j\downarrow0$. Let $f_j$ be almost full ball charts
centered at $p-\delta_jv$, with inverses
$G_j\to G:\BB^n\to\partial\Om$ locally uniformly
and $G(0)=p$. Assume $G$ has constant rank $r$ near $0$,
so its image germ at $p$ is embedded.
Suppose also that
\[
 F_j(s)=f_j(p-\delta_jsv)\longrightarrow F(s),\qquad
 F:\HH\to\BB^n,\qquad F(1)=0,
\]
locally uniformly on $\HH$, with the maps $F_j$ eventually defined
on every compact subset. Choose a sufficiently small holomorphic
product neighborhood $U$ of $p$, and let $M$ be the zero section
representing the image germ in $U$. There is a compact
piece $M_0\Subset M\cap U$, containing $p$, and a number
$\theta\in(0,1)$ with the following property.
For every neighborhood $W\Subset U$ of $M_0$, put
$V=\Om\cap W$. Then, for every $Q\Subset\HH$ and all
sufficiently large $j$,
\begin{equation}\label{eq:fixed-leaf-localization}
 1\le\frac{K_V(p-\delta_jsv)}{K_\Om(p-\delta_jsv)}
 \le\theta^{-2n}\qquad(s\in Q).
\end{equation}
The same $\theta$ works for arbitrarily thin $W$.
\end{lem}

\begin{proof}
\proofstep{1}{A fixed compact subset of the leaf} To make the bound independent of thickness, we use
balls of fixed radius whose images stay in a compact part of the leaf.
For fixed $s\in\HH$, $F_j(s)$ stays in a compact subset of $\BB^n$,
so $G_j(F_j(s))=p-\delta_j s v$ gives $G(F(s))=p$.
Let $\mathcal W$ be the component of $G^{-1}(U)$ containing $0$.
The defining transverse coordinates of $M$ vanish on $G$ near $0$,
hence $G(\mathcal W)\subset M$ by the identity theorem. The connected set $F(\HH)$ contains
$0$ and maps to $p$, so it lies in $\mathcal W$.

Choose $R>0$ with $G(\BB^n)\subset B(p,R)$ and
$0<\theta<\theta'<1$ with $\overline{B(p,R\theta')}\Subset U$.
Fix a compact zero-section piece $M_0\Subset M\cap U$ such that
\[
 \overline{B(p,R\theta')}\cap M\subset\operatorname{int}_M M_0.
\]
Fix these before $Q$ and $W$.
For $a\in\BB^n$, let $A_a$ be a ball automorphism sending
$a$ to $0$. For $a=F(s)$ we have $G(a)=p$, so the ball Schwarz lemma applied to
\[
 z\longmapsto\frac{G(A_a^{-1}(z))-p}{R}
\]
gives
\[
 G\bigl(A_{F(s)}^{-1}(\theta'\BB^n)\bigr)\subset B(p,R\theta')\subset U.
\]
The source ball is connected and contains $F(s)\in\mathcal W$,
so it lies in $\mathcal W$. Its image thus lies in $B(p,R\theta')\cap M\subset M_0$.

\proofstep{2}{Ball inclusions and kernel bounds} For $s\in Q$, the closed balls
$A_{F(s)}^{-1}(\theta'\overline{\BB^n})$ lie in the closed
Kobayashi ball about $0$ of radius
$\sup_Q k_{\BB^n}(0,F(s))+\operatorname{arctanh}\theta'$, by the triangle
inequality. Uniform convergence gives
$k_{\BB^n}(F_j(s),F(s))<\operatorname{arctanh}\theta'-
\operatorname{arctanh}\theta$ on $Q$ for large $j$, so
\[
 A_{F_j(s)}^{-1}(\theta\BB^n)
 \subset A_{F(s)}^{-1}(\theta'\BB^n)
 \qquad(s\in Q).
\]
For large $j$, the common compact lies in $r_j\BB^n$.
There $G_j\to G$ uniformly, and $M_0$ has positive distance
from the complement of $W$, giving
\[
 G_j\bigl(A_{F_j(s)}^{-1}(\theta\BB^n)\bigr)\subset V
 \qquad(s\in Q)
\]
for large $j$. With $D_j=f_j(\Om)$ and $a_j=F_j(s)$,
\eqref{eq:kernel-monotone}, \eqref{eq:kernel-transform},
and \eqref{eq:ball-kernel} give
\[
 K_{f_j(V)}(a_j)
 \le K_{A_{a_j}^{-1}(\theta\BB^n)}(a_j)
 =\theta^{-2n}K_{\BB^n}(a_j)
 \le\theta^{-2n}K_{D_j}(a_j).
\]
Transforming back gives the upper bound in
\eqref{eq:fixed-leaf-localization}; inclusion $V\subset\Om$
gives the lower bound. Use componentwise kernels on disconnected open sets.
\end{proof}

The following translation estimate is common to the comparison in
a thin tube and the comparison after $G_j\to p$.
Only the passage from local to global kernels differs.

\begin{lem}[Translation in nested neighborhoods]
\label{lem:nested-translation}
In a holomorphic coordinate neighborhood of $p$, suppose
$\Psi(0,w)=p+w\nu_p$ and the local defining inequality is
$\Rea w+\phi(\zeta,\Ima w)<0$, with $\zeta$ denoting all
the tangential coordinates. Let
$\widehat W_0\Subset\widehat W_1$ be product neighborhoods of
the coordinate origin and assume $|\partial_t\phi|\le\varepsilon$
on $\widehat W_1$, where $t=\Ima w$.
Set $V_i=\Om\cap\Psi(\widehat W_i)$. For every compact
$Q\subset\{x+iy:x>\varepsilon|y|\}$, uniformly for $s=x+iy$ in $Q$,
\begin{align}
 K_{V_0}(q_\delta(s))&\ge(1-o(1))
       K_{V_1}(q_\delta(x+\varepsilon|y|)),\\
 K_{V_0}(q_\delta(x-\varepsilon|y|))&\ge(1-o(1))
       K_{V_1}(q_\delta(s))
\end{align}
as $\delta\downarrow0$.
\end{lem}

\begin{proof}
Write $g(\zeta,w)=\Rea w+\phi(\zeta,\Ima w)$,
$\xi_\delta(s)=(0,-\delta s)$ and $\widehat V_i=\Psi^{-1}(V_i)$.
The three points tend uniformly to the coordinate origin on $Q$.
Since $\phi(0,0)=0$, the real points have $g<0$ by $x\pm\varepsilon|y|>0$.
At $\xi_\delta(s)$, $g(\xi_\delta(s))\le-\delta(x-\varepsilon|y|)<0$.
Thus all three lie in $\widehat V_0$ for small $\delta$.
Define
\[
 T_{\delta,\pm y}(\zeta,w)
   =(\zeta,w-\varepsilon\delta|y|\pm i\delta y).
\]
The displacement is $\delta|y|\sqrt{1+\varepsilon^2}$, less than the
distance from $\overline{\widehat W_0}$ to the complement of $\widehat W_1$
for small $\delta$, uniformly on $Q$. Along the segment,
$|\partial_t\phi|\le\varepsilon$
gives $g\circ T_{\delta,\pm y}-g\le-\varepsilon\delta|y|+\varepsilon\delta|y|=0$.
Thus $T_{\delta,\pm y}(\widehat V_0)\subset\widehat V_1$.
At the relevant normal points,
\begin{align*}
 T_{\delta,y}(\xi_\delta(s))&=\xi_\delta(x+\varepsilon|y|),\\
 T_{\delta,-y}(\xi_\delta(x-\varepsilon|y|))&=\xi_\delta(s).
\end{align*}
Translations have determinant one, so domain monotonicity gives the
two inequalities in the coordinate neighborhoods without error.
In the original coordinates, divide by the Jacobian factors
$|\det d\Psi(\xi_\delta(s'))|^2$. Their ratios tend uniformly to one
since each factor tends to $|\det d\Psi(0)|^2>0$.
\end{proof}

\begin{lem}[Vertical comparison with a fixed constant]
\label{lem:tube-comparison}
Under the assumptions of \Cref{lem:local-leaf-localization},
suppose that $r>0$, $v=\nu_p$, $\Om$ is pseudoconvex, and its
boundary is $C^2$ near $p$. There is $C_0\ge1$ such that,
for every $\varepsilon>0$,
\begin{equation}\label{eq:local-fixed-comparison}
 C_0^{-1}(1-o(1))K_\Om(q_{\delta_j}(x+\varepsilon|y|))
 \le K_\Om(q_{\delta_j}(x+iy))
 \le C_0(1+o(1))K_\Om(q_{\delta_j}(x-\varepsilon|y|))
\end{equation}
uniformly on compact subsets of $\{x>\varepsilon|y|\}$.
The constant $C_0$ depends on the fixed charts, their limit,
and the initial tube, not on $\varepsilon$ or the compact set; the threshold for $j$ may depend on both.
\end{lem}

\begin{proof}
Fix the base piece $M_0$ and $\theta$ in the tube supplied by
\Cref{lem:tube,lem:local-leaf-localization}. To make
$|\partial_t\phi|<\varepsilon$, shrink only the fibers and choose
nested product tubes
$\widehat W_0\Subset\widehat W_1$ whose zero sections contain
$M_0$. Put $W_i=\Psi(\widehat W_i)$ and $V_i=\Om\cap W_i$.
The localization lemma gives
\[
 1\le\frac{K_{V_i}(q_{\delta_j}(s))}
              {K_\Om(q_{\delta_j}(s))}\le C_0,
 \qquad C_0=\theta^{-2n},\quad i=0,1,
\]
on every fixed compact subset of the half-plane for all large $j$.
Further shrinking the fibers leaves $C_0$ unchanged.

Apply \Cref{lem:nested-translation}. Its first inequality gives
the lower bound using $C_0^{-1}K_{V_0}\le K_\Om$ at $s$ and
$K_\Om\le K_{V_1}$ at $x+\varepsilon|y|$. The second gives the upper
bound using $K_\Om\le K_{V_1}$ at $s$ and $K_{V_0}\le C_0K_\Om$
at $x-\varepsilon|y|$, proving \eqref{eq:local-fixed-comparison}.
\end{proof}

The chosen subsequence of scales and charts stays fixed. Shrinking the tube fibers changes only the index
from which the estimates hold, not $\theta$ or
$C_0=\theta^{-2n}$. Thus one may first let $j\to\infty$
and then $\varepsilon\downarrow0$ on that same subsequence.

\subsection{Boundary curves and convergence of inverse charts}\label{sec:exclude-curves}

If a boundary curve exists, the thin-tube comparison carries
the real-axis lower bound for the kernel ratios to compact
subsets of $\HH$. The tube also supplies holomorphic discs
with a nonzero derivative in the leaf parameter. The next
lemma shows that this derivative is incompatible with the
resulting lower bound for the normalized kernel limit. Once
boundary curves are excluded, every inverse-chart limit is constant.

\begin{lem}[A superharmonic obstruction]\label{lem:product-obstruction}
Let $\mathcal H:\HH\times\DD\to\BB^n$ be holomorphic,
put $F(s)=\mathcal H(s,0)$, and suppose
$A(s)=\partial_\zeta\mathcal H(s,0)$ is nowhere zero.
Let $J$ be a nowhere-zero holomorphic function on $\HH$ and set
\[
 k(s)=\frac{|J(s)|^2}{(1-\norm{F(s)}^2)^{n+1}}.
\]
There are no constants $\mu>0$ and $C\in\R$ such that
\begin{equation}\label{eq:forbidden-profile-blowup}
 \log k(x+iy)\ge\mu\log(1/x)-C
 \qquad(0<x\le1,\ |y|\le1).
\end{equation}
\end{lem}

\begin{proof}
On $|\zeta|=r<1$, the powers of $\zeta$ are orthogonal.
Keeping only the constant and linear terms gives
\[
 \norm{F(s)}^2+r^2\norm{A(s)}^2
 \le \frac{1}{2\pi}\int_0^{2\pi}
       \norm{\mathcal H(s,re^{i\theta})}^2\,d\theta\le1.
\]
Let $r\uparrow1$, so $\norm{A}^2\le1-\norm{F}^2$.
Since $J$ has no zeros, $\log|J|^2$ is harmonic.
Subtracting the subharmonic $(n+1)\log\norm{A}^2$ gives the
smooth superharmonic function $U=\log|J|^2-(n+1)\log\norm{A}^2\ge\log k$.
If \eqref{eq:forbidden-profile-blowup} held, then
$U+C\ge0$ on $0<x\le1$, $|y|\le1$. On the left side of
$Q_\epsilon=\{x+iy:\epsilon<x<1,\ |y|<1\}$ it would
be at least $\mu\log(1/\epsilon)$.
The fixed harmonic function
\[
 b(x,y)=\frac{\sinh(\frac\pi2(1-x))}{\sinh(\pi/2)}\cos(\tfrac\pi2 y)
\]
vanishes on the other three sides and satisfies $0\le b\le1$ on the
left side. By the minimum principle,
$U+C-\mu\log(1/\epsilon)b\ge0$ on $Q_\epsilon$. For $0<\epsilon<1/2$, this gives at $s=1/2$
\[
 U(\tfrac12)+C\ge\mu\log(1/\epsilon)b(\tfrac12,0)\longrightarrow+\infty,
 \qquad b(\tfrac12,0)>0,
\]
a contradiction.
\end{proof}

\begin{thm}[Exclusion of boundary curves]\label{thm:no-curves}
Let $\Om\Subset\C^n$ be a bounded pseudoconvex domain, and assume
only that its boundary is $C^2$ in a neighborhood of
$p\in\partial\Om$. If \eqref{eq:normal-hypothesis} holds,
then no nonconstant holomorphic curve germ through $p$ is
contained in $\partial\Om$.
\end{thm}

\begin{proof}
\proofstep{1}{Select normal scales} Suppose a boundary curve germ
exists. The rough bounds in \Cref{lem:kernel-bounds} give
\[
 2t-C\le h(t):=\log K_\Om(q_{e^{-t}}(1))
 \le(n+1)t+C.
\]
Use \Cref{lem:envelope} with $b=2$, $\mu=1$, and $B>n+1$.
The full normal limit supplies almost full ball charts at $\delta_j=e^{-T_j}$.
Pass to subsequences so that $G_j\to G$ and $F_j\to F$,
where $F_j(s)=f_j(q_{\delta_j}(s))$. The maximum principle
gives $F:\HH\to\BB^n$ with $F(1)=0$.

\Cref{lem:curve-rank}
gives $r:=\rank dG(0)>0$, and \Cref{prop:local-rank} shows that
$G$ has constant rank near $0$. Its image germ $M$ is therefore
a nonsingular complex manifold in the boundary near $p$. \Cref{lem:tube-comparison} applies with a
constant $C_0$ independent of the fiber thickness.

\proofstep{2}{Pass to the limit on $\HH$} Fix $Q\Subset\HH$, choose $\varepsilon>0$
so that $\inf_Q(x-\varepsilon|y|)>0$. Then $x+\varepsilon|y|$
ranges over a compact subset of $(0,\infty)$, where the two real-axis
bounds in \eqref{eq:asymmetric-profile-bound} are uniformly positive.
The lower bound in \eqref{eq:local-fixed-comparison} therefore bounds
$R_j(s)=K_\Om(q_{\delta_j}(s))/K_\Om(q_{\delta_j}(1))$
away from zero on $Q$ for all sufficiently large $j$. Hence
\Cref{lem:jacobian-compactness} gives a nowhere-zero
holomorphic $J$ and the kernel limit
\[
 R_j(s)\longrightarrow
 k(s)=\frac{|J(s)|^2}{(1-\norm{F(s)}^2)^{n+1}}.
\]
In the lower bound of \eqref{eq:local-fixed-comparison}, let
$j\to\infty$ and then $\varepsilon\downarrow0$, keeping $C_0$ fixed:
\begin{equation}\label{eq:vertical-comparability}
 k(x+iy)\ge C_0^{-1}k(x)
 \qquad(x>0,\ y\in\R).
\end{equation}
Since \eqref{eq:asymmetric-profile-bound} gives $k(x)\ge x^{-1}$
for $0<x\le1$, we obtain
\begin{equation}\label{eq:uniform-profile-blowup}
 \log k(x+iy)\ge\log(1/x)-\log C_0
 \qquad(0<x\le1,\ y\in\R).
\end{equation}

\proofstep{3}{The superharmonic obstruction} Use the same tube and choose a small nonzero
$u\in\C^r$ so that $\gamma(\zeta)=g(\zeta u)$ is defined
on a neighborhood of $\overline\DD$. Since $g$ is an immersion,
$\gamma'(0)\ne0$. On compact subsets of $\HH\times\DD$,
the maps
\[
 \mathcal H_j(s,\zeta)
   =f_j\bigl(\Psi(\zeta u,0,-\delta_js)\bigr)
\]
are eventually defined. Indeed,
$\phi(\eta,0,0)=\partial_t\phi(\eta,0,0)=0$ gives the
defining inequality $-\delta_j\Rea s+O(\delta_j^2)<0$
on each such compact set. After passing to a subsequence, they converge to
$\mathcal H:\HH\times\DD\to\BB^n$; the value
$\mathcal H(1,0)=0$ ensures that the entire limit image lies in the open ball.
Setting $\zeta=0$ and using \eqref{eq:tube-normal} gives the
first identity below. Applying $G_j$ to $\mathcal H_j$ and passing
to the limit gives the second.
\[
 \mathcal H(s,0)=F(s),\qquad
 G\circ\mathcal H(s,\zeta)=\gamma(\zeta).
\]
Thus $A(s)=\partial_\zeta\mathcal H(s,0)$ is nowhere zero,
because $dG(F(s))A(s)=\gamma'(0)\ne0$.
Now \eqref{eq:uniform-profile-blowup} contradicts
\Cref{lem:product-obstruction}. This excludes the curve germ.
\end{proof}

\begin{cor}\label{cor:inverse-collapse}
Under the domain and local $C^2$ assumptions of
\Cref{thm:no-curves} and \eqref{eq:normal-hypothesis}, let
$\delta_j\downarrow0$
and let $f_j:\Om\hookrightarrow\BB^n$ satisfy
$f_j(q_{\delta_j}(1))=0$ and
$r_j\BB^n\subset f_j(\Om)$, with $r_j\to1$.
Then the inverses $G_j=f_j^{-1}$ converge locally uniformly
on $\BB^n$ to the constant map $p$.
\end{cor}

\begin{proof}
Boundedness gives a normal family.
By \Cref{lem:boundary-limit,thm:no-curves}, every subsequential
limit is constant and takes $0$ to $p$. Thus $G_j\to p$.
\end{proof}

\section{Normalized kernel limits and translation semigroups}\label{sec:profiles}

Under local $C^2$ regularity, \Cref{thm:no-curves} excludes boundary
curves when \eqref{eq:normal-hypothesis} holds. We find a normalized
kernel limit whose holomorphic factor $J$ satisfies
$J(1+t)=\det d\Phi_t(0)$ for the limiting semigroup and $t>0$.

\subsection{Localization using the same normal charts}

Once $G_j\to p$, the same charts can be used to compare the
kernels on $\Om$ and $\Om\cap U$ for a sufficiently small fixed
neighborhood $U$ of $p$. Along the subsequence on which the
normal-line maps converge, the ratio tends to one uniformly
on compact subsets of the rescaled half-plane. No
plurisubharmonic peak function is needed. As in \cite{FNlocal}, such
a function for a domain $D\subset\C^m$ at $q\in\partial D$ is a function
$\varphi\in\operatorname{PSH}(D)$ such that $\varphi(z)\to0$ as
$D\ni z\to q$ and $\sup_{D\setminus U}\varphi<0$ for every
neighborhood $U$ of $q$. At such a peak point $p$, Nikolov and Verma
\cite[Proposition~2]{NV} show that $s_\Om(z)\to1$ at $p$ implies
$s_{\Om\cap U}(z)\to1$ at $p$ for bounded $U$ with $\Om\cap U$
connected. Forn\ae ss and Nikolov restate this in
\cite[Theorem~1(ii), p.~355]{FNlocal}; their remark explains why the
converse needs global hypotheses. The proof below does not use these
squeezing-localization results.

\begin{lem}\label{lem:normal-localization}
Let $\Om\Subset\C^n$ be a bounded domain with $C^2$ boundary
near $p\in\partial\Om$, and suppose that no boundary curve germ
through $p$ exists. Let
$f_j$ be almost full ball charts centered at $q_{\delta_j}(1)$,
where $\delta_j\downarrow0$, and put $D_j=f_j(\Om)$.
Then the inverses $G_j=f_j^{-1}$ converge locally uniformly to
$p$ on $\BB^n$. After passing to a subsequence,
\[
 F_j(s)=f_j(q_{\delta_j}(s))\longrightarrow F(s)
 \quad\text{locally uniformly on }\HH,\qquad F(1)=0,
\]
where $F:\HH\to\BB^n$ is holomorphic. For every sufficiently
small fixed neighborhood $U$ of $p$,
\begin{equation}\label{eq:normal-localization}
 \frac{K_{\Om\cap U}(q_{\delta_j}(s))}
      {K_\Om(q_{\delta_j}(s))}\longrightarrow1
 \qquad\text{locally uniformly on }\HH.
\end{equation}
\end{lem}

\begin{proof}
Boundedness of $\Om$ gives subsequential limits of $G_j$ on $\BB^n$.
By \Cref{lem:boundary-limit} each limit maps into $\partial\Om$ and
takes $0$ to $p$. A nonconstant limit would yield a boundary curve germ
through $p$ by restriction to a suitable complex line through $0$.
Thus $G_j\to p$ along the whole sequence.

The normal-line maps $F_j$ take values in the unit ball and are
eventually defined on each compact subset of $\HH$.
Montel gives a locally uniform limit on a subsequence. Since $F(1)=0$,
the maximum principle puts its entire image in the open ball.

Fix $Q\Subset\HH$ and $\sup_Q\norm{F}<r<1$. Then $F_j(Q)\subset r\BB^n$
eventually. For each $r<a<1$, convergence of $G_j$ on
$a\overline{\BB^n}$ gives $a\BB^n\subset E_j\subset D_j\subset\BB^n$
eventually, where $E_j=f_j(\Om\cap U)$.
Both kernels transform under $f_j$ with the same Jacobian
factor, which cancels in the quotient.
Thus for $s\in Q$,
\[
 1\le
 \frac{K_{\Om\cap U}(q_{\delta_j}(s))}
      {K_\Om(q_{\delta_j}(s))}
 =
 \frac{K_{E_j}(F_j(s))}{K_{D_j}(F_j(s))}
 \le
 \frac{K_{a\BB^n}(F_j(s))}{K_{\BB^n}(F_j(s))}.
\]
The final ratio is
$a^{-2n}[(1-|u|^2)/(1-|u|^2/a^2)]^{n+1}$ at $u=F_j(s)$
by \eqref{eq:kernel-transform} and \eqref{eq:ball-kernel}.
Take the supremum over $s\in Q$ and then $\limsup_{j\to\infty}$.
Since the right-hand ratio tends uniformly to one on $|u|\le r$
as $a\uparrow1$, the kernel quotient converges uniformly to one on $Q$.
Only the charts centered at $q_{\delta_j}(1)$ are used, not squeezing
at nonreal parameters.
\end{proof}

\begin{lem}\label{lem:normal-comparison}
In the setting of \Cref{lem:normal-localization}, for every
$\varepsilon>0$ one has
\begin{equation}\label{eq:normal-comparison}
 (1-o(1))K_\Om(q_{\delta_j}(x+\varepsilon|y|))
 \le K_\Om(q_{\delta_j}(x+iy))
 \le(1+o(1))K_\Om(q_{\delta_j}(x-\varepsilon|y|)),
\end{equation}
uniformly on compact subsets of $\{x>\varepsilon|y|\}$.
\end{lem}

\begin{proof}
First compare locally. In affine normal coordinates the graph is $\Rea w+\phi(z,t)<0$, where $t=\Ima w$ and
$\phi(0)=d\phi(0)=0$.
Choose nested product neighborhoods $U_0\Subset U_1$ on which
$|\partial_t\phi|<\varepsilon$. Apply \Cref{lem:nested-translation} with $\Psi=\id$. By \Cref{lem:normal-localization},
the local-to-global kernel ratios at the three parameters
$x+iy$ and $x\pm\varepsilon|y|$ tend to one uniformly on the
relevant compact sets. Replacing the local kernels in the
two inequalities gives \eqref{eq:normal-comparison}.
\end{proof}

\subsection{Selecting normal scales from kernel growth}

The comparison transfers the real-axis kernel bounds to
each compact subset of $\HH$, giving the positive lower bound
needed for compactness of the reciprocal Jacobians. The choice
of $B$ depends on the application.

\begin{prop}[A normalized kernel limit along selected scales]\label{prop:selected-profile}
Let $\Om\Subset\C^n$ be a bounded pseudoconvex domain with $C^2$
boundary near $p$ and no boundary curve germ through $p$. Assume
\eqref{eq:normal-hypothesis}, and put
\[
 h(t)=\log K_\Om(q_{e^{-t}}(1)).
\]
Let $B>2$ satisfy
\begin{equation}\label{eq:profile-envelope-condition}
 \liminf_{t\to\infty}\bigl(h(t)-Bt\bigr)=-\infty.
\end{equation}
Then there exist scales $\delta_j\downarrow0$ and almost
full ball charts $f_j$ whose inverses tend locally uniformly
to $p$ and for which the maps $F_j(s)=f_j(q_{\delta_j}(s))$
and the normalized Jacobians $J_j$ from \eqref{eq:jacobian-ratio}
have the following limits locally uniformly on $\HH$:
\[
 F_j\longrightarrow F:\HH\longrightarrow\BB^n,
 \qquad J_j\longrightarrow J,
 \qquad F(1)=0,\quad J(1)=1,
\]
where $J$ is holomorphic and nowhere zero. The normalized
kernels have the $C^\infty_{\mathrm{loc}}(\HH)$ limit
\begin{equation}\label{eq:normal-profile}
 k(s)=\frac{|J(s)|^2}{(1-|F(s)|^2)^{n+1}},
\end{equation}
and
\begin{equation}\label{eq:profile-bounds}
 k(x)\ge
 \begin{cases}
 x^{-1},&0<x\le1,\\
 x^{-B},&x\ge1.
 \end{cases}
\end{equation}
\end{prop}

\begin{proof}
We select scales on the real axis, then extend the lower bound to
complex parameters. By \Cref{lem:kernel-bounds}, $h(t)\ge2t-C$.
Apply \Cref{lem:envelope} with $b=2$ and the stated $B$, and use
\eqref{eq:normal-hypothesis} to choose charts at $\delta_j=e^{-T_j}$.
The ratios satisfy \eqref{eq:asymmetric-profile-bound} on the real axis;
\Cref{lem:normal-localization} gives $G_j\to p$ and
subsequential convergence $F_j\to F$ on $\HH$.

For $Q\Subset\HH$ choose $\varepsilon>0$ with $x-\varepsilon|y|>0$ on $Q$.
Then $x+\varepsilon|y|$ ranges in a compact subinterval of $(0,\infty)$,
where the lower bound in \eqref{eq:asymmetric-profile-bound} is positive
uniformly. The lower inequality in \eqref{eq:normal-comparison} gives
$R_j\ge c_Q>0$ on $Q$ for large $j$. Now \Cref{lem:jacobian-compactness}
gives $J$ and \eqref{eq:normal-profile}; passing to the limit in
\eqref{eq:asymmetric-profile-bound} gives \eqref{eq:profile-bounds}.
\end{proof}

\subsection{Normal translations on the full ball}\label{sec:translation}

A holomorphic semigroup on a complex manifold $X$ is a family
$(\Phi_t)_{t\ge0}$ of holomorphic self-maps of $X$, continuous in
the compact-open topology, with $\Phi_0=\id_X$ and
$\Phi_t\circ\Phi_u=\Phi_{t+u}$ for nonnegative real times.
Its generator is the holomorphic vector field whose flow gives the semigroup.

To recover the full $n$-dimensional Jacobian, we take limits
of the conjugated normal translations on the whole ball, not
only on the complex normal line. The factor $J$ in
\eqref{eq:normal-profile} then
satisfies $J(1+t)=\det d\Phi_t(0)$ for $\Rea t>0$.

\begin{prop}\label{prop:translation-semigroup}
For the scales and charts in \Cref{prop:selected-profile}, a
further subsequence of
\begin{equation}\label{eq:Phi-j}
 \Phi_j(t,z)=f_j\bigl(G_j(z)-\delta_jt\nu_p\bigr)
\end{equation}
converges locally uniformly on $\HH\times\BB^n$ to a
holomorphic map $\Phi$. With $\Phi_t(z)=\Phi(t,z)$, it
satisfies
\begin{align}
 \Phi_t\circ\Phi_u&=\Phi_{t+u}
       &&(t,u\in\HH),\label{eq:semigroup-law}\\
 \Phi_t(0)&=F(1+t),\qquad
 \det d\Phi_t(0)=J(1+t)\ne0.\label{eq:exact-Jacobian}
\end{align}
Moreover, $t\mapsto\Phi_t$ extends continuously to the
closed right half-plane in the compact-open topology. In particular,
\begin{equation}\label{eq:Phi-zero}
 \Phi_t\longrightarrow\id_{\BB^n}
 \qquad(t\to0,\ \Rea t>0).
\end{equation}
For every $y\in\R$,
\begin{equation}\label{eq:imaginary-limit}
 A_y=\lim_{\varepsilon\downarrow0}\Phi_{\varepsilon+iy}
\end{equation}
exists in that topology, and $(A_y)_{y\in\R}$ is a continuous
one-parameter subgroup of $\Aut(\BB^n)$.
\end{prop}

\begin{proof}
\proofstep{1}{Translations inside the domain}
To work on the whole ball in the limit, fix $K\Subset\BB^n$ and
$T\Subset\HH$. Since $G_j(K)\to p$, a local defining function satisfies
\begin{equation}\label{eq:Phi-domain}
 \rho(G_j(z)-\delta_jt\nu_p)
 =\rho(G_j(z))-2\delta_j\Rea\bigl(t\partial\rho_{G_j(z)}(\nu_p)\bigr)
 +O_{K,T}(\delta_j^2).
\end{equation}
Write $a=\partial\rho_p(\nu_p)>0$ and $m=\inf_{t\in T}\Rea t>0$
for nonempty $T$. Since $T$ is bounded and the conormal pairing
converges uniformly to $a$ on $K$, the real part in the linear term
is at least $am/2$ for all large $j$. Thus the inward term of size
at least $am\delta_j$ dominates the quadratic remainder. This proves
eventual definition on each product compact set without a rate for
$G_j\to p$.

\proofstep{2}{The limit and its Jacobian}
Exhaust $\HH\times\BB^n$ by relatively compact open products.
On each product the maps are defined and bounded by one for all large $j$.
Montel's theorem and diagonal extraction give a jointly holomorphic
limit. At the base point the exact identities are
\[
 \Phi_j(t,0)=F_j(1+t),\qquad
 \det d_z\Phi_j(t,0)
 =\det df_j(q_{\delta_j}(1+t))\det dG_j(0)=J_j(1+t).
\]
Convergence with derivatives yields
\eqref{eq:exact-Jacobian}. Since $F(1+t)$ lies in the open
ball, the maximum principle shows that $\Phi_t$ maps the
entire ball into the open ball. In particular, each $\Phi_t$
is locally biholomorphic at $0$, so its image contains a
neighborhood of $\Phi_t(0)$.

\proofstep{3}{The composition law}
For $t,u\in\HH$ and $K\Subset\BB^n$, the sets $\Phi_j(u,K)$
eventually lie in one compact subset of $\BB^n$, hence in $r_j\BB^n$.
The inward-translation estimate allows composition there, and
$G_j\circ f_j=\id$ gives
\[
 \Phi_j(t,\Phi_j(u,z))
 =f_j(G_j(z)-\delta_j(t+u)\nu_p)=\Phi_j(t+u,z)
 \qquad(z\in K).
\]
Pass to the limit on that compact set to obtain \eqref{eq:semigroup-law}.

\proofstep{4}{Extension to imaginary time}
Fix a real $u>0$. If $t_l\in\HH$ tends to $t$ in the closed half-plane,
Montel gives a compact-open cluster limit $A$ of $\Phi_{t_l}$.
Every such limit is ball-valued,
since $A(0)=F(1+t)\in\BB^n$. By \eqref{eq:semigroup-law},
$A\circ\Phi_u=\Phi_{t+u}$. The image of $\Phi_u$ contains an open set,
so the identity theorem determines $A$ uniquely. This gives the
boundary value $A_y$ at $t=iy$ and $A_0=\id$.

The same identity and base-point value now hold at every closed time.
For any sequence of closed times approaching $t$, they force every
cluster limit to equal the value at $t$. Thus the extension is continuous.
Compact-open convergence allows composition: the inner maps send each
compact set into a common compact subset of the ball for large indices.
Taking interior-time limits extends \eqref{eq:semigroup-law} to closed
times. In particular, $A_yA_v=A_{y+v}$ and
$A_yA_{-y}=A_{-y}A_y=A_0=\id$.
\end{proof}

\begin{remark}\label{rem:complex-time-order}
The finite-$j$ maps are only required for $\Rea t>0$.
No assertion is made that a purely imaginary physical
translation preserves $\Om$. The imaginary automorphisms
arise after first taking the chart limit in the open
half-plane, and then taking its boundary values.
\end{remark}

\section{Semigroup rigidity and normal kernel growth}\label{sec:classification}

\subsection{Classification and kernel ratios}

For a holomorphic ball self-map without an interior fixed point,
its iterates converge locally uniformly in the ball to the same
boundary point, called its Denjoy--Wolff point; see Herv\'e's ball theorem
\cite[Theorem~4, p.~3-9]{Her63}.
We use the normal form for a single automorphism in
Bracci, Contreras, and D\'iaz-Madrigal
\cite[Proposition~4.3]{BCDclass}, then derive the complex-time
restrictions and the decay of a normalized kernel ratio.

\begin{thm}\label{thm:semigroup-classification}
Let $\Phi:\HH\times\BB^n\to\BB^n$ be holomorphic, with
$\Phi_t\Phi_u=\Phi_{t+u}$ for $t,u\in\HH$. In the compact-open
topology, assume $\Phi_t\to\Phi_0:=\id$ as $t\to0$ in $\HH$ and that
the boundary values
\[
 A_y=\lim_{\varepsilon\downarrow0}\Phi_{\varepsilon+iy}
 \qquad(y\in\R)
\]
exist and form a continuous one-parameter group in $\Aut(\BB^n)$.
After a fixed biholomorphic conjugacy, its real-time maps
$\Phi_t$, $t\ge0$, have one of the following forms:
\begin{enumerate}[label=\textup{(\roman*)},leftmargin=*]
\item on $\BB^n$, with $D$ an $n\times n$ Hermitian matrix,
\begin{equation}\label{eq:elliptic-form}
 \Phi_t(z)=e^{-tD}z,\qquad D=D^*\ge0;
\end{equation}
\item on the Siegel domain
\begin{equation}\label{eq:Siegel-domain}
 \mathcal S=\{(z,w)\in\C^{n-1}\times\C:
                    \Rea w>|z|^2\},
\end{equation}
with $D$ an $(n-1)\times(n-1)$ Hermitian matrix,
\begin{equation}\label{eq:Siegel-form}
 \Phi_t(z,w)=(e^{-tD}z,w+ct),\qquad D=D^*\ge0,\quad c>0.
\end{equation}
\end{enumerate}
At a fixed $a\in\BB^n$, define the kernel ratio including the Jacobian factor:
\begin{equation}\label{eq:density}
 \mathcal V_a(t)=
 |\det d\Phi_t(a)|^2
 \frac{K_{\BB^n}(\Phi_t(a))}{K_{\BB^n}(a)}.
\end{equation}
If the semigroup is not the identity semigroup, then
\begin{equation}\label{eq:density-decay}
 \mathcal V_a(t)=O(t^{-(n+1)})\qquad(t\to+\infty).
\end{equation}
More precisely, it decays at least exponentially when $D\ne0$;
when $D=0$ in \eqref{eq:Siegel-form},
it is exactly a positive constant multiple of the power function
$t^{-(n+1)}$ after a shift of its argument.
\end{thm}

\begin{proof}
\proofstep{1}{The real-time generator}
Let $X$ be the holomorphic generator of the continuous automorphism
group. For $t>0$ the semigroup law and the boundary limit give
$A_y\Phi_t=\lim_{\varepsilon\downarrow0}\Phi_{t+\varepsilon+iy}
=\Phi_{t+iy}$. Differentiating in $y$ at $0$, for
each fixed real $t>0$, yields
\[
 X\circ\Phi_t=i\partial_t\Phi_t.
\]
Thus the real-time flow solves the equation for $V=-iX$.
Integrate along an orbit from $\varepsilon$ to $t$ and use
$\Phi_\varepsilon\to\id$ as $\varepsilon\downarrow0$ to obtain the
initial condition $\Phi_0=\id$, without differentiating at the time boundary.

\proofstep{2}{A common interior fixed point}
To put $V$ in normal form, apply the common-fixed-point alternative
\cite[Theorem~3.1, p.~35]{BCD2010} to the semigroup $(A_y)_{y\ge0}$.
If it has a common fixed point in the ball, conjugate that point
to $0$. Every $A_y$ is then unitary, so
\[
 A_y(z)=e^{iyH}z,\qquad H=H^*.
\]
The vector field $V=-iX$ is $V(z)=Hz$, and uniqueness of its
flow gives $\Phi_t(z)=e^{tH}z$. Preservation of the ball for
all $t\ge0$ forces every eigenvalue of $H$ to be nonpositive.
Set $D=-H$ to obtain \eqref{eq:elliptic-form}.

\proofstep{3}{A common boundary fixed point}
Otherwise the same theorem gives a common boundary Denjoy--Wolff
point for $A_y$ with $y>0$; the inverses also fix it.
Rotate it to $e_n$, then send it to infinity in
\eqref{eq:Siegel-domain} by the Cayley map $\mathcal C$, whose inverse is
\begin{equation}\label{eq:inverse-cayley}
 \mathcal C^{-1}(z,w)=\left(\frac{2z}{w+1},\frac{w-1}{w+1}\right).
\end{equation}
Each $A_y$, $y>0$, now has no interior fixed point and has its
Denjoy--Wolff point at infinity. Since ball automorphisms are
linear fractional, \cite[Proposition~4.3, formula~(8), p.~333]{BCDclass}
gives the following expression for $A_y$ for $y>0$:
\[
 (z,w)\longmapsto
 \bigl(r_yU_yz+q_y,\,
 r_y^2w+2r_y\ip{U_yz}{q_y}+|q_y|^2+i\beta_y\bigr),
\]
where $r_y>0$, $U_y$ is unitary, $q_y\in\C^{n-1}$, and
$\beta_y\in\R$. The continuous one-parameter subgroup is smooth
in the automorphism group. The displayed parameters are
recovered from the values and derivative blocks at $(0,1)$,
so they are differentiable from the right at $0$, with the
identity values $(r_0,U_0,q_0,\beta_0)=(1,I,0,0)$.
Differentiating the formula with respect to $y$ at zero gives
\begin{equation}\label{eq:Siegel-X}
 X(z,w)=\bigl((\ell I+iH)z+a_0,
              2\ell w+2\ip{z}{a_0}+ib\bigr),
\end{equation}
where $\ell=r'_0$, $iH=U'_0$ with $H=H^*$, $a_0=q'_0$,
and $b=\beta'_0$. Recall that the inner product is linear in its first variable.
Consequently
\[
 V(z,w)=-iX(z,w)
 =\bigl(Hz-i\ell z-ia_0,
         -2i\ell w-2i\ip{z}{a_0}+b\bigr).
\]
For $h(z,w)=\Rea w-|z|^2$, use
$dh(V)=\Rea V_w-2\Rea\ip{V_z}{z}$ on $\partial\mathcal S$:
\begin{equation}\label{eq:inward-polynomial}
 dh(V)=b+2\ell\Ima w-2\ip{Hz}{z}+4\Ima\ip{z}{a_0}\ge0.
\end{equation}
Indeed, the affine field extends its local flow across each Euclidean
boundary point. The interior flow preserves $\mathcal S$, so continuity
gives $h\ge0$ along its boundary trajectories for small positive times,
and the initial right derivative is nonnegative.

On $\partial\mathcal S$, the variable $\Ima w$ is arbitrary.
Thus \eqref{eq:inward-polynomial} forces $\ell=0$.
The quadratic part in $z$ then gives $D:=-H\ge0$.
For $u\in\ker D$, substitute all complex multiples of $u$;
nonnegativity forces $\ip{u}{a_0}=0$. Hence
$a_0\in(\ker D)^\perp=\ran D$. The Moore--Penrose inverse
$D^\dagger$ equals $(D|_{\ran D})^{-1}$ on $\ran D$ and zero
on $\ker D$. Set $c_0=iD^\dagger a_0$.
Completing the square gives
\begin{equation}\label{eq:completion-square}
 b+2\ip{Dz}{z}+4\Ima\ip{z}{a_0}
 =2\ip{D(z+c_0)}{z+c_0}+b-2\ip{D^\dagger a_0}{a_0}.
\end{equation}
It follows that $c:=b-2\ip{D^\dagger a_0}{a_0}\ge0$.
The Heisenberg coordinate change
\begin{equation}\label{eq:Heisenberg-shift}
 Z=z+c_0,\qquad
 W=w+2\ip{z}{c_0}+|c_0|^2
\end{equation}
preserves $h$. Differentiating along $V$ gives
\[
 VZ=-DZ,\qquad VW=c.
\]
The flow has the form \eqref{eq:Siegel-form}. Since $X=iV$,
if $c=0$ the shifted group $A_y$ would be
$(Z,W)\mapsto(e^{-iyD}Z,W)$, fixing every $(0,W)$ with
$\Rea W>0$. This contradicts this case of the
fixed-point alternative. Hence $c>0$. The two normal forms
are distinguished by whether the imaginary-time group has
a common interior fixed point.

\proofstep{4}{The kernel ratio in each normal form}
By \eqref{eq:kernel-transform}, the Jacobian factors of the conjugacy
at the base point and its image cancel in \eqref{eq:density}.
It therefore suffices to compute in the displayed coordinates.

In \eqref{eq:elliptic-form}, write $v$ for the conjugated
base point. Insert $\det e^{-tD}=e^{-t\tr D}$ and
\eqref{eq:ball-kernel} into \eqref{eq:density}:
\begin{equation}\label{eq:elliptic-density}
 \mathcal V_a(t)=e^{-2t\tr D}
 \left(\frac{1-|v|^2}{1-|e^{-tD}v|^2}\right)^{n+1}
 \le e^{-2t\tr D}.
\end{equation}
If $D\ne0$, then $\tr D>0$; otherwise the semigroup is the
identity.

In \eqref{eq:Siegel-form}, write the conjugated point as
$(v,w_0)$ and put $h_0=\Rea w_0-|v|^2>0$.
For the Siegel kernel, the map \eqref{eq:inverse-cayley} gives
\[
 1-|\mathcal C^{-1}(z,w)|^2=\frac{4h(z,w)}{|w+1|^2},
 \qquad
 \det d\mathcal C^{-1}(z,w)=\frac{2^n}{(w+1)^{n+1}}.
\]
Equations \eqref{eq:kernel-transform} and \eqref{eq:ball-kernel} give
\begin{equation}\label{eq:Siegel-kernel}
 K_{\mathcal S}(z,w)=\frac{n!}{4\pi^n}\,h(z,w)^{-(n+1)}.
\end{equation}
It follows that
\begin{equation}\label{eq:Siegel-density}
 \mathcal V_a(t)=e^{-2t\tr D}
 \left(
 \frac{h_0}{\Rea w_0+ct-|e^{-tD}v|^2}
 \right)^{n+1}.
\end{equation}
The denominator is at least $h_0+ct$. If $D\ne0$, we again
obtain exponential decay. If $D=0$, the formula reduces
to $(1+ct/h_0)^{-(n+1)}$, with $c>0$.
This proves \eqref{eq:density-decay}.
\end{proof}

\begin{remark}\label{rem:fixed-base}
The decay in \eqref{eq:density-decay} need not be uniform in $a$.
For a pure Siegel translation and each fixed $t>0$,
$(1+ct/h_0)^{-(n+1)}\to1$ as $h_0\to\infty$. We use only $a=0$.
\end{remark}

\begin{prop}[The exact normalized Bergman kernel limit]\label{prop:exact-selected-profile}
Every normalized kernel limit $k$ from \Cref{prop:selected-profile} satisfies
\begin{equation}\label{eq:exact-selected-profile}
 k(s)=(\Rea s)^{-(n+1)}\qquad(s\in\HH).
\end{equation}
The corresponding limits $F,J$ have the form below for a unit vector $\xi\in\C^n$:
\begin{equation}\label{eq:exact-orbit-Jacobian}
 F(s)=\frac{s-1}{s+1}\xi,\qquad
 J(s)=\left(\frac{2}{s+1}\right)^{n+1}
 \qquad(s\in\HH).
\end{equation}
\end{prop}

\begin{proof}
\proofstep{1}{Identify the normal form} Put $N=n+1$. By \Cref{prop:translation-semigroup} and
\eqref{eq:normal-profile}, \eqref{eq:exact-Jacobian},
\begin{equation}\label{eq:density-identity}
 k(1+t)=|\det d\Phi_t(0)|^2
 \frac{K_{\BB^n}(\Phi_t(0))}{K_{\BB^n}(0)}
 =\mathcal V_0(t)\qquad(t>0).
\end{equation}
As $t\to\infty$, the lower bound $k(1+t)\ge(1+t)^{-B}$ in
\eqref{eq:profile-bounds} rules out exponential decay in
\Cref{thm:semigroup-classification}. The identity semigroup
would give $F\equiv0$ and $J\equiv1$ by
\eqref{eq:exact-Jacobian} and the holomorphic identity theorem,
contrary to the small-$x$ lower bound. Thus there is a fixed
biholomorphism $T:\BB^n\to\mathcal S$ such that
\[
 T\circ\Phi_t\circ T^{-1}(z,w)=(z,w+ct),\qquad c>0.
\]
Write $T(0)=(v,w_0)$ and $h_0=\Rea w_0-|v|^2>0$.
Evaluate the conjugacy and its Jacobian at $0$, with $s=1+t>1$
and \eqref{eq:exact-Jacobian}. Translation has determinant one,
so the identity theorem extends both holomorphic identities to $\HH$:
\begin{equation}\label{eq:conjugated-normal-orbit}
 T(F(s))=(v,w_0+c(s-1)),\qquad
 J(s)\,\frac{\det dT(F(s))}{\det dT(0)}=1.
\end{equation}
The Siegel height of this orbit is $h_0+c(\Rea s-1)>0$.
Letting real $s\downarrow0$ gives $h_0\ge c$.
The Bergman transformation
formula and \eqref{eq:Siegel-kernel} now yield
\begin{equation}\label{eq:conjugated-profile-density}
 k(s)=\left(\frac{h_0}{h_0+c(\Rea s-1)}\right)^N
 \qquad(s\in\HH).
\end{equation}
If $h_0>c$, the right side remains bounded as real
$s=x\downarrow0$, contrary to $k(x)\ge x^{-1}$.
Thus $h_0=c$, proving \eqref{eq:exact-selected-profile}.

\proofstep{2}{The orbit and Jacobian} The Heisenberg translation
$P(z,w)=(z-v,w-2\langle z,v\rangle+|v|^2-i\Ima w_0)$
and dilation $L(z,w)=(c^{-1/2}z,c^{-1}w)$ send
$T(0)$ to $(0,1)$ and $T(F(s))$ to $(0,s)$.
Composing $LPT$ with the inverse Cayley map
\eqref{eq:inverse-cayley} gives a ball automorphism $U$
that fixes zero and is therefore unitary.
Thus $UF(s)=(s-1)(s+1)^{-1}e_n$. The Jacobian determinants
of $U$ and $L\circ P$ are constant, so the determinant identity above gives
\[
 J(s)=\frac{\det d\mathcal C^{-1}(0,s)}{\det d\mathcal C^{-1}(0,1)}
      =\left(\frac{2}{s+1}\right)^N
\]
by \eqref{eq:inverse-cayley}, proving \eqref{eq:exact-orbit-Jacobian}.
\end{proof}

Projecting the ball onto the disc $F(\HH)$ gives a holomorphic
retraction. Hence $F$ preserves the Kobayashi distance and
infinitesimal metric. The formulas concern only selected subsequences.
\Cref{lem:normal-comparison} extends the real-axis lower bounds
to ensure $J$ is nowhere zero in \Cref{prop:selected-profile}.

\subsection{Normal kernel growth and the main theorem}\label{sec:conclusion}

The formula for $k$ forces the logarithmic growth exponent to
equal $n+1$. At a Levi-degenerate point, the order-$2+\alpha$
remainder bound \eqref{eq:pointwise-taylor} would give a strictly
smaller exponent.

\begin{proof}[Proof of \Cref{thm:normal-log-growth}]
Write $N=n+1$ and $h(t)=\log K_\Om(q_{e^{-t}}(1))$.
The kernel bounds give
\begin{equation}\label{eq:log-growth-rough}
 2t-C\le h(t)\le Nt+C.
\end{equation}
If $L:=\liminf_{t\to\infty}h(t)/t<N$, choose $\max\{2,L\}<B<N$.
Along a sequence $t_j\to\infty$ with $h(t_j)/t_j\to L$,
eventually $h(t_j)-Bt_j\le-(B-L)t_j/2\to-\infty$.
Thus \eqref{eq:profile-envelope-condition} holds. With no boundary
curve germ by \Cref{thm:no-curves}, \Cref{prop:selected-profile} gives $k$.
By \Cref{prop:exact-selected-profile}, $k(2)=2^{-N}$,
whereas \eqref{eq:profile-bounds} gives $k(2)\ge2^{-B}$.
Since $B<N$, these inequalities contradict each other.
Thus $\liminf h(t)/t\ge N$, and the upper bound in
\eqref{eq:log-growth-rough} proves
\eqref{eq:normal-log-growth}.
\end{proof}

\begin{proof}[Proof of \Cref{thm:main}]
If the Levi form has nullity $\ell\ge1$, \Cref{lem:kernel-bounds}
gives $K_\Om(q_\delta(1))\le C\delta^{-a_{\alpha,\ell}}$, where
$a_{\alpha,\ell}=n+1-\ell\alpha/(2+\alpha)
\le n+2/(2+\alpha)<n+1$.
Taking logs and dividing by $\log(1/\delta)$ gives an upper limit
at most $a_{\alpha,\ell}$, since $\log C/\log(1/\delta)\to0$.
This contradicts \Cref{thm:normal-log-growth}.
\end{proof}

\begin{proof}[Proof of \Cref{cor:local-question}]
At each relevant point the boundary limit gives
\eqref{eq:normal-hypothesis}, so \Cref{thm:main} applies.
\end{proof}

\begin{cor}[A pointwise Taylor condition]\label{cor:pointwise-taylor}
Let $\Om\Subset\C^n$ be a bounded pseudoconvex domain, with
$C^2$ boundary near $p\in\partial\Om$. Suppose that
some $C^2$ defining function near $p$ satisfies
\eqref{eq:pointwise-taylor} for some $0<\alpha\le1$.
Then \eqref{eq:normal-hypothesis} implies that $p$ is
strictly pseudoconvex.
\end{cor}

\begin{proof}
At a Levi-degenerate point, the remainder bound
\eqref{eq:pointwise-taylor} gives the same kernel upper bound
\eqref{eq:weak-kernel-bound}. The logarithmic comparison
in the proof of \Cref{thm:main} again contradicts \Cref{thm:normal-log-growth}.
\end{proof}

The remainder bound \eqref{eq:pointwise-taylor} is imposed at
$p$ alone; the second derivatives need not be H\"older continuous
on a neighborhood of $p$.

\subsection{Limit order and the regularity threshold}

In \Cref{thm:no-curves}, take the chart limit on fixed compact
sets before $x\downarrow0$ in the superharmonic obstruction;
this does not use a semigroup. After $G_j\to p$,
\Cref{prop:translation-semigroup} identifies $J(1+t)$ with
$\det d\Phi_t(0)$ for $t>0$. For the selected limit $k$,
\Cref{prop:exact-selected-profile} combines the semigroup classification
with the lower bounds as $x\to\infty$ and $x\downarrow0$.
\Cref{thm:normal-log-growth} then uses the formula at $x=2$.

The conormal, curve exclusion, and logarithmic kernel growth
use only local $C^2$ regularity. The order-$2+\alpha$ bound
\eqref{eq:pointwise-taylor} on the Taylor remainder of a defining
function is used in \Cref{thm:main,cor:pointwise-taylor} only to
obtain the improved kernel upper bound \eqref{eq:weak-kernel-bound}
at a Levi-degenerate point.

At the $C^2$ threshold, if the Levi form at $p$ is degenerate
and \eqref{eq:normal-hypothesis} holds, then
\Cref{lem:kernel-bounds,thm:normal-log-growth} give
the simultaneous limits
\begin{equation}\label{eq:critical-kernel-behavior}
 \delta^{n+1}K_\Om(p-\delta\nu_p)\longrightarrow0,
 \qquad
 \frac{\log K_\Om(p-\delta\nu_p)}{\log(1/\delta)}
 \longrightarrow n+1.
\end{equation}
For every $\varepsilon>0$, the logarithmic limit also gives
$\delta^{n+1}K_\Om(p-\delta\nu_p)\ge\delta^\varepsilon$
at sufficiently small depths. Thus the rescaled kernel tends to zero,
but more slowly than any positive power of $\delta$.
The function $\delta^{-(n+1)}/\log(1/\delta)$ illustrates these
growth rates; no domain with this exact kernel is asserted.
Thus $C^2$ regularity alone does not give a kernel upper exponent
strictly below $n+1$ at a Levi-degenerate point. This is consistent
with the $C^2$ counterexample
of Forn\ae ss--Wold \cite{FW}.

\section{Sampling and geometric applications}\label{sec:applications}

We derive sampling criteria and realization obstructions from
the main theorem and the preceding $C^2$ results. In particular,
under the domain and local boundary-regularity assumptions of
\Cref{thm:normal-log-growth},
\[
 \liminf_{\delta\downarrow0}
 \frac{\log K_\Om(p-\delta\nu_p)}{\log(1/\delta)}<n+1
 \quad\Longrightarrow\quad
 \liminf_{\delta\downarrow0}s_\Om(p-\delta\nu_p)<1.
\]
Indeed, if the conclusion failed, then
\eqref{eq:normal-hypothesis} would hold, contradicting that theorem.
Thus any geometric estimate giving a strict gap in the normal
kernel growth exponent yields an obstruction to squeezing tending to one.

\begin{cor}[Anisotropic polydisc obstruction]
\label{cor:polydisc-obstruction}
Let $\Om\Subset\C^n$, $n\ge2$, be a bounded pseudoconvex
domain with $C^2$ boundary near $p$. In affine unitary
coordinates with $p=0$ and $\nu_p$ the last coordinate vector,
suppose that along a sequence $\delta_m\downarrow0$,
\[
 \bigl\{(z,w): |z_j|<c\delta_m^{\beta_j}\ (1\le j\le n-1),
              \ |w+\delta_m|<c\delta_m\bigr\}\subset\Om
\]
for fixed $c>0$ and $\beta_1,\ldots,\beta_{n-1}>0$.
If
\[
 a:=2+2\sum_{j=1}^{n-1}\beta_j<n+1,
\]
then $\liminf_{\delta\downarrow0}s_\Om(p-\delta\nu_p)<1$.
\end{cor}

\begin{proof}
To force a growth gap, compare with the inscribed polydisc,
whose center kernel is the reciprocal of its volume. By monotonicity,
\[
 K_\Om(q_{\delta_m}(1))\le\pi^{-n}c^{-2n}\delta_m^{-a}.
\]
The logarithmic quotient has liminf at most $a<n+1$.
The squeezing obstruction follows from
\Cref{thm:normal-log-growth}.
\end{proof}

For example, tangential radii $\delta_m^{1/(2m_j)}$ with
$m_j\ge1$ and at least one $m_j>1$ give
$a=2+\sum_jm_j^{-1}<n+1$. The conclusion concerns the liminf along the full
normal; the argument does not identify a squeezing
gap along the prescribed sequence $\delta_m$.

The same obstruction holds if the displayed polydiscs $P_m$
are replaced by $\Psi(P_m)$, where $\Psi$ is a fixed local
biholomorphism near $0$ satisfying $\Psi(0,w)=(0,w)$.
Indeed, the transformation law adds the factor
$|\det d\Psi(q_{\delta_m}(1))|^{-2}$ to the kernel bound,
which is bounded above and below by positive constants.
This includes the shear $(z,w)\mapsto(z,w+Q(z))$ for a
homogeneous holomorphic quadratic polynomial $Q$.

An interior tangent ball gives an upper bound for the
infinitesimal Kobayashi metric in the normal direction and
controls the distance between successive normal samples.

\begin{lem}[Comparison along an interior tangent ball]
\label{lem:normal-ball-comparison}
Suppose $B=B(p-R\nu_p,R)\subset\Om$ is an interior tangent
ball. For $0<h<R$, write $p_h=p-h\nu_p$ and use the
logarithmic diameter coordinate $L_R(h)=\log(h/(2R-h))$. Then $d_\Om(p_h)=h$, and
\begin{align}
 h\kappa_\Om(p_h;\nu_p)&\le\frac{1}{2-h/R},
 \label{eq:normal-ball-metric}\\
 e^{2 k_\Om(p_h,p_k)}
 &\le e^{|L_R(h)-L_R(k)|}
 \le 2\max\{h/k,k/h\}\qquad(0<h,k<R).
 \label{eq:normal-ball-distance}
\end{align}
\end{lem}

\begin{proof}
We compute on the tangent ball and then use inclusion.
The radius-$h$ ball about $p_h$ lies in $B$, while
$p\in\partial\Om$ is at distance $h$. Scale $B$ to the unit
ball and project orthogonally onto the complex normal disc;
this is a holomorphic left inverse of the disc inclusion.
The coordinate of $p_h$ is $1-h/R$ and normal speed is $1/R$.
Use $\kappa_\DD(x;v)=|v|/(1-x^2)$ and
$k_\DD(x,y)=|\operatorname{arctanh}x-\operatorname{arctanh}y|$
for $x,y\in(-1,1)$ to obtain
\[
 \kappa_B(p_h;\nu_p)=\frac{1}{h(2-h/R)},\qquad
 2k_B(p_h,p_k)=|L_R(h)-L_R(k)|.
\]
The inclusion $B\subset\Om$ decreases both metrics and distances,
so it gives \eqref{eq:normal-ball-metric} and $k_\Om\le k_B$ on $B$.
For the distance bound, when $h\le k$, the exponential of the
last expression is $k(2R-h)/(h(2R-k))\le 2k/h$;
interchanging $h,k$ handles the other order.
\end{proof}

The formula for the selected limit map $F$ also gives the
limits of the Kobayashi distance and of the infinitesimal
Kobayashi and Bergman metrics along the complex normal.

\begin{cor}[Intrinsic geometry along selected scales]
\label{cor:selected-intrinsic-geometry}
Let $\Om\Subset\C^n$, $n\ge2$, be a bounded pseudoconvex
domain with $C^2$ boundary near $p\in\partial\Om$, and assume
\eqref{eq:normal-hypothesis}. There is a sequence
$\delta_j\downarrow0$ such that
\begin{align}
 k_\Om(q_{\delta_j}(s),q_{\delta_j}(t))
 &\longrightarrow k_\HH(s,t),
 \label{eq:selected-normal-distance}\\
 \delta_j\kappa_\Om(q_{\delta_j}(s);\nu_p)
 &\longrightarrow\frac{1}{2\Rea s},
 \label{eq:selected-normal-kobayashi}\\
 \delta_j b_\Om(q_{\delta_j}(s);\nu_p)
 &\longrightarrow\frac{\sqrt{n+1}}{2\Rea s}.
 \label{eq:selected-normal-bergman}
\end{align}
The first convergence is locally uniform on $\HH\times\HH$,
and the other two are locally uniform on $\HH$. Moreover,
\begin{equation}\label{eq:normal-kobayashi-limsup}
 \limsup_{\delta\downarrow0}\delta\kappa_\Om(q_\delta(1);\nu_p)=\frac12,
 \qquad
 \limsup_{\delta\downarrow0}\delta b_\Om(q_\delta(1);\nu_p)
 =\frac{\sqrt{n+1}}2.
\end{equation}
\end{cor}

\begin{proof}
\proofstep{1}{Compare the Kobayashi and Bergman metrics}
Put $N=n+1$ and compare the metrics on almost full ball images $r_j\BB^n\subset D_j\subset\BB^n$,
$r_j\to1$. On each fixed compact subset of $\BB^n$, the bounds
$\kappa_{\BB^n}(a;v)\le\kappa_{D_j}(a;v)\le
\kappa_{\BB^n}(a/r_j;v/r_j)$ converge uniformly for $|v|=1$.
Since $K_{D_j}\ge K_{\BB^n}>0$, \Cref{lem:ball-kernel} gives uniform
convergence of the complex Hessians of $\log K_{D_j}$ to that of
$\log K_{\BB^n}$ on these compact sets. At $0$, $b_{\BB^n}(0;v)^2=N|v|^2$ by \eqref{eq:ball-kernel},
and $\kappa_{\BB^n}(0;v)=|v|$ by the Schwarz lemma.
Ball automorphisms give $b_{\BB^n}=\sqrt N\,\kappa_{\BB^n}$.
The ball Kobayashi metric is bounded below by a positive
constant on unit vectors. Thus, by homogeneity,
\[
 \frac{b_{D_j}(a;v)}{\kappa_{D_j}(a;v)}
 \longrightarrow\sqrt N
\]
uniformly for $a$ in each such compact set and all $v\ne0$.

\proofstep{2}{Limits along selected scales} Use $B>N$: the estimate
$h(t)=\log K_\Om(q_{e^{-t}}(1))\le Nt+C$ from
\Cref{lem:kernel-bounds} gives $h(t)-Bt\to-\infty$.
With \Cref{thm:no-curves}, this permits the selection in
\Cref{prop:selected-profile,prop:exact-selected-profile}, giving
\[
 F_j(s)=f_j(q_{\delta_j}(s))\longrightarrow
 F(s)=\frac{s-1}{s+1}\xi,\qquad |\xi|=1.
\]
Write $D_j=f_j(\Om)$, with
$r_j\BB^n\subset D_j\subset\BB^n$ and $r_j\to1$.
For $Q\Subset\HH$, eventually $F_j(Q)\subset r_j\BB^n$.
Biholomorphic invariance and domain monotonicity give, for $s,t\in Q$,
\[
 k_{\BB^n}(F_j(s),F_j(t))
 \le k_\Om(q_{\delta_j}(s),q_{\delta_j}(t))
 \le k_{\BB^n}(F_j(s)/r_j,F_j(t)/r_j).
\]
Both bounds converge uniformly to $k_\HH(s,t)$, since $F$
is a complex geodesic. This proves \eqref{eq:selected-normal-distance}.

For the infinitesimal bounds, Cauchy's estimates give
$F'_j\to F'$ locally uniformly and the chain rule gives
$F'_j(s)=-\delta_jdf_j(q_{\delta_j}(s))\nu_p$. Thus
\[
 \kappa_{\BB^n}(F_j(s);F'_j(s))
 \le\delta_j\kappa_\Om(q_{\delta_j}(s);\nu_p)
 \le\kappa_{\BB^n}(F_j(s)/r_j;F'_j(s)/r_j).
\]
Both bounds converge locally uniformly to
$\kappa_\HH(s;1)=(2\Rea s)^{-1}$.
The vectors $F'_j(s)$ are nonzero because $f_j$ is biholomorphic.
By invariance and homogeneity,
$b_\Om(q_{\delta_j}(s);\nu_p)/\kappa_\Om(q_{\delta_j}(s);\nu_p)
=b_{D_j}(F_j(s);F'_j(s))/\kappa_{D_j}(F_j(s);F'_j(s))$.
On each $Q\Subset\HH$, the points $F_j(Q)$ remain in a fixed
compact subset of $\BB^n$. Combining the ratio comparison
with \eqref{eq:selected-normal-kobayashi} proves
\eqref{eq:selected-normal-bergman}.

\proofstep{3}{Upper limits along the normal} Local $C^2$ regularity gives
an interior tangent ball. By \eqref{eq:normal-ball-metric},
$\limsup_{\delta\downarrow0}\delta\kappa_\Om(q_\delta(1);\nu_p)\le1/2$.
Along the selected sequence, \eqref{eq:selected-normal-kobayashi} at
$s=1$ gives the limit $1/2$, so the upper limit equals $1/2$.
For the metric ratio, now take any positive sequence $\delta_j\to0$.
The squeezing limit along the entire normal gives almost full
ball charts centered at $q_{\delta_j}(1)$ with radii $r_j\to1$.
Applying Step~1 at the chart center $a=0$ and using biholomorphic
invariance give
$b_\Om(q_{\delta_j}(1);\nu_p)/\kappa_\Om(q_{\delta_j}(1);\nu_p)
\to\sqrt N$.
By the tangent-ball bound, the normalized Kobayashi metric stays bounded,
so $\delta(b_\Om-\sqrt N\,\kappa_\Om)(q_\delta(1);\nu_p)=o(1)$
along the full normal. Taking upper limits proves the
second identity in \eqref{eq:normal-kobayashi-limsup}.
For a quantitative comparison at an arbitrary interior
point in terms of its squeezing function, see Nikolov and
Trybu\l a \cite[Section~II, before Corollary~3]{NT20}.
\end{proof}

At strictly pseudoconvex points, full boundary asymptotics for
these metrics already hold under $C^2$ regularity; see
\cite[Theorems~2 and~4 and the remarks following them]{NT20}.
The corollary gives the Kobayashi distance and rescaled metric
limits along selected scales without assuming strict pseudoconvexity.
Without the squeezing hypothesis, even bounded smooth pseudoconvex
domains can have a sequence with normalized normal Kobayashi metric
tending to zero \cite[Theorem~3]{FL09}.

\subsection{Sampling the normal and changing approach curves}
\label{subsec:normal-sampling}

Ball automorphisms bound the squeezing defect at one point by
the defect at another times an exponential in their Kobayashi
distance. The bound controls the errors even when the distances
between successive samples grow. For related continuity estimates,
see \cite[Theorems~3.1--3.2]{DGZbasic}.

\begin{lem}[Comparing squeezing at two points]\label{lem:squeezing-transport}
Let $D\subset\C^n$ be a bounded domain, let $z,w\in D$, and put
$\tau=\tanh k_D(z,w)$. Then
\begin{equation}\label{eq:squeezing-transport}
 \max\left\{0,\frac{s_D(z)-\tau}{1-\tau s_D(z)}\right\}
 \le s_D(w)
 \le \frac{s_D(z)+\tau}{1+\tau s_D(z)}.
\end{equation}
In particular,
\begin{equation}\label{eq:squeezing-defect-transport}
 1-s_D(w)\le e^{2k_D(z,w)}\bigl(1-s_D(z)\bigr).
\end{equation}
Consequently, if $s_D(z_j)\to1$ and
$\sup_j k_D(z_j,w_j)<\infty$, then $s_D(w_j)\to1$.
\end{lem}

\begin{proof}
Set $d=k_D(z,w)<\infty$ and $S(x)=\operatorname{arctanh}s_D(x)$,
with $S(x)=+\infty$ when $s_D(x)=1$. We first prove $S(w)\ge S(z)-d$.
For $0<r<s_D(z)$,
choose an embedding $f$ centered at $z$ whose image contains
$r\BB^n$, and put $a=f(w)$. Then $k_{\BB^n}(0,a)\le d$.
Since $k_{\BB^n}(0,\xi)=\operatorname{arctanh}|\xi|$, if
$\operatorname{arctanh}r>d$, the triangle inequality
shows that the ball of Kobayashi radius
$\operatorname{arctanh}r-d$ about $a$ lies in $r\BB^n$.
Moving $a$ to $0$ by an automorphism gives
$S(w)\ge\operatorname{arctanh}r-d$. If $\operatorname{arctanh}r\le d$,
the same bound follows from $S(w)\ge0$. Let $r\uparrow s_D(z)$,
then exchange $z,w$ to obtain
\[
 S(z)-d\le S(w)\le S(z)+d
\]
in the extended real order. Since $d$ is finite, either both
values of $S$ are infinite or both are finite. Apply the addition
formula for $\tanh$, extended continuously at $+\infty$, and use
$s_D(w)\ge0$ to obtain \eqref{eq:squeezing-transport}.

Subtract the untruncated lower bound in \eqref{eq:squeezing-transport}
from one. Since $1-\tau s_D(z)\ge1-\tau>0$, this gives
\[
 1-s_D(w)\le
 \frac{(1+\tau)(1-s_D(z))}{1-\tau s_D(z)}
 \le e^{2d}(1-s_D(z)),
 \qquad \frac{1+\tau}{1-\tau}=e^{2d}.
\]
Bounded $d$ then preserves convergence to one.
\end{proof}

\begin{lem}[A nontangential point and its normal comparison point]
\label{lem:nontangential-normal-distance}
Let $D\subset\C^n$ be a domain with $C^1$ boundary near $p\in\partial D$,
and let $\nu_p$ be the outward unit normal. For every $A\ge1$
there are $C_A<\infty$ and a neighborhood $U_A$ of $p$ such
that, whenever $z\in D\cap U_A$ satisfies
\begin{equation}\label{eq:nontangential-cone}
 |z-p|\le A d_D(z),
\end{equation}
the point $p-d_D(z)\nu_p$ belongs to $D$ and
\begin{equation}\label{eq:cone-normal-distance}
 k_D\bigl(z,p-d_D(z)\nu_p\bigr)\le C_A.
\end{equation}
\end{lem}

\begin{proof}
We find a segment whose distance to $\partial D$ is at least
$d_D(z)/4$. In real orthogonal coordinates $(x,u)$ with
$p=0$ and $\nu_p$ the positive $u$-direction,
\[
 D=\{(x,u):u+\phi(x)<0\},\qquad
 \phi(0)=0,\quad d\phi(0)=0.
\]
Choose $0<\eta\le\min\{1,1/(4A)\}$ and shrink the coordinate
neighborhood so that $|d\phi|\le\eta$ there.
Put $\delta=d_D(z)$ and write $z=(x,u)$. For $z$ near zero, $(x,-\phi(x))$ lies on the local boundary. Thus
\[
 \delta\le-u-\phi(x),\qquad
 -u\ge\delta-\eta|x|\ge(1-\eta A)\delta.
\]
For the segment
\[
 \gamma(t)=((1-t)x,(1-t)u-t\delta),\qquad0\le t\le1,
\]
we therefore have
\[
 -\gamma_u(t)-\phi(\gamma_x(t))
 \ge(1-2\eta A)\delta\ge\delta/2.
\]
The function $(x,u)\mapsto u+\phi(x)$ is Lipschitz with
constant at most $\sqrt{1+\eta^2}$ on a fixed ball about the
origin in the coordinate neighborhood. Since $\delta\le |z|$, every
ball of radius $\delta/4$ centered on the segment lies within
distance $5|z|/4$ of the origin, hence in that fixed ball for small $|z|$.
The defining function changes there by at most
$\sqrt{1+\eta^2}\,\delta/4<\delta/2$ from the segment center, so
each such ball lies in $D$.
Join successive endpoints by holomorphic discs. If $\ell=|z+\delta\nu_p|=0$,
the distance is zero. Otherwise take $m$ equal pieces with
$4\ell/(m\delta)<1$.
In the original complex coordinates, put
$v=(-\delta\nu_p-z)/\ell$. On the piece starting at $a$,
the affine holomorphic disc $\zeta\mapsto a+\zeta(\delta/4)v$
joins the endpoints at parameters $0$ and $4\ell/(m\delta)$.
These discs lie in $D$ by the ball containment just proved.
The sum of the Poincar\'e distances bounds the Kobayashi distance:
\[
 k_D(z,-\delta\nu_p)
 \le m\operatorname{arctanh}\!\left(\frac{4\ell}{m\delta}\right)
 \longrightarrow \frac{4\ell}{\delta}\le4(A+1)
 \qquad(m\to\infty).
\]
Thus $C_A=4(A+1)$ works after shrinking $U_A$.
\end{proof}

\begin{cor}[Normal convergence from successive samples]
\label{cor:normal-sampling}
Let $\Om\Subset\C^n$, $n\ge2$, be a bounded pseudoconvex domain with
$C^{2,\alpha}$ boundary near $p\in\partial\Om$ for some
$0<\alpha\le1$.
Suppose $z_j\to p$ nontangentially, and set
$\delta_j=d_\Om(z_j)$ and $\varepsilon_j=1-s_\Om(z_j)$.
Assume that $\delta_j\downarrow0$ and
\begin{equation}\label{eq:quantitative-normal-sampling}
 \varepsilon_j\varepsilon_{j+1}
 \frac{\delta_j}{\delta_{j+1}}\longrightarrow0.
\end{equation}
Then $s_\Om(p-\delta\nu_p)\to1$ as $\delta\downarrow0$,
and $p$ is strictly pseudoconvex. In particular, the condition
holds if $s_\Om(z_j)\to1$ and
$\delta_{j+1}\ge c\delta_j$ for some $c>0$.
For normal samples $z_j=p-h_j\nu_p$ with $h_j\downarrow0$,
the boundary distance $\delta_j$ equals $h_j$ for all large $j$.
\end{cor}

\begin{proof}
We bound the squeezing defect at each intermediate normal
point using the defects at both sample points. Fix $A\ge1$ with $|z_j-p|\le A\delta_j$
eventually; \Cref{lem:nontangential-normal-distance} gives
$k_\Om(z_j,p-\delta_j\nu_p)\le C_A$.
Fix an interior tangent ball of radius $R$ at $p$.
By \Cref{lem:normal-ball-comparison}, normal points at depth
$h<R$ have boundary distance $h$, and $L=L_R$ is increasing. For $\delta_{j+1}\le\delta\le\delta_j<R$,
the triangle inequality through $p-\delta_k\nu_p$ and
\eqref{eq:normal-ball-distance} give
$2k_\Om(z_k,p-\delta\nu_p)\le2C_A+|L(\delta_k)-L(\delta)|$.
By \eqref{eq:squeezing-defect-transport},
\[
 1-s_\Om(p-\delta\nu_p)
 \le e^{2C_A}\varepsilon_k e^{|L(\delta_k)-L(\delta)|},
 \qquad k\in\{j,j+1\}.
\]
Multiply the two bounds: the exponents add to
$L(\delta_j)-L(\delta_{j+1})$, so their exponential product
is at most $2\delta_j/\delta_{j+1}$. Taking square roots gives
\begin{equation}\label{eq:sampling-uniform-defect}
 \sup_{\delta\in[\delta_{j+1},\delta_j]}
 \bigl(1-s_\Om(p-\delta\nu_p)\bigr)
 \le\sqrt2\,e^{2C_A}
 \sqrt{\varepsilon_j\varepsilon_{j+1}
              \frac{\delta_j}{\delta_{j+1}}}
 \longrightarrow0.
\end{equation}
The squeezing function therefore tends to one along the entire
normal, and \Cref{thm:main} applies.
\end{proof}

\begin{remark}\label{rem:sampling-scope}
The product condition \eqref{eq:quantitative-normal-sampling}
already forces $\varepsilon_j\to0$ in this geometric setting.
Apply the endpoint bound with $\delta=\delta_j$ and $k=j+1$.
The distance from $p-\delta_j\nu_p$ to $z_j$ is at most $C_A$,
so transferring the defect to $z_j$ by
\eqref{eq:squeezing-defect-transport} adds a factor $e^{2C_A}$.
Since $e^{L(\delta_j)-L(\delta_{j+1})}\le2\delta_j/\delta_{j+1}$,
\[
 \varepsilon_j^2\le2e^{4C_A}
 \varepsilon_j\varepsilon_{j+1}\frac{\delta_j}{\delta_{j+1}}
 \longrightarrow0.
\]
The distance control is essential; an abstract condition on products
of adjacent terms alone would not imply convergence of the factors.

A sufficient special case is $s_\Om(z_j)\to1$ together with
bounded ratios of successive depths, or equivalently bounded successive gaps
in $-\log\delta_j$. This covers the choices
$\delta_j=e^{-j}$ and $\delta_j=j^{-1}$. Unlike the sequence criteria at $h$-extendible points
\cite{Nik,MV,NN20}, this result requires no finite-type assumption,
but it relates sampling gaps to squeezing defects; it does not
cover arbitrary sparse sequences.

Since $0\le\varepsilon_{j+1}\le1$, the one-endpoint condition
$\varepsilon_j\delta_j/\delta_{j+1}\to0$ implies the two-endpoint
condition, but not conversely. For example, $\varepsilon_j=1/j$ and
$\delta_j/\delta_{j+1}=j$ for $j\ge2$ make the one-endpoint
product equal to one, whereas the product in
\eqref{eq:quantitative-normal-sampling} equals $1/(j+1)$.
This compares numerical conditions, not their realization by a domain.
\end{remark}

\begin{cor}[Continuous nontangential approaches]
\label{cor:nontangential-curve}
Let $\Om\Subset\C^n$, $n\ge2$, be a bounded pseudoconvex domain with
$C^{2,\alpha}$ boundary near $p\in\partial\Om$ for some
$0<\alpha\le1$. Suppose
$\gamma:(0,t_0]\to\Om$ is continuous,
$\gamma(t)\to p$ as $t\downarrow0$, and
\[
 |\gamma(t)-p|\le A d_\Om(\gamma(t))
\]
for all sufficiently small $t$ and some $A<\infty$.
If $s_\Om(\gamma(t))\to1$, then $p$ is strictly
pseudoconvex.
In particular, if $\gamma:[0,t_0]\to\overline\Om$ is $C^1$,
$\gamma(0)=p$, $\gamma(t)\in\Om$ for $t>0$,
$s_\Om(\gamma(t))\to1$ as $t\downarrow0$, and
\[
 \Rea\ip{\gamma'(0)}{\nu_p}<0,
\]
then $p$ is strictly pseudoconvex.
\end{cor}

\begin{proof}
We compare points at equal distances from the boundary. Choose $t_1>0$ so
the nontangential bound and \Cref{lem:nontangential-normal-distance}
apply on $(0,t_1]$. The continuous function $d_\Om\circ\gamma$
tends to zero at $0$. For each $0<\delta<d_\Om(\gamma(t_1))$,
the intermediate value theorem gives $t(\delta)\in(0,t_1)$
with $d_\Om(\gamma(t(\delta)))=\delta$.
Every choice has $t(\delta)\to0$ as $\delta\downarrow0$, since
$d_\Om\circ\gamma$ has a positive minimum on each $[a,t_1]$, $0<a\le t_1$.
The equal-depth comparison and \eqref{eq:squeezing-defect-transport} give
\[
 1-s_\Om(p-\delta\nu_p)
 \le e^{2C_A}\bigl(1-s_\Om(\gamma(t(\delta)))\bigr)\longrightarrow0.
\]
Apply \Cref{thm:main} at $p$.

For the final assertion, put $a=-\Rea\ip{\gamma'(0)}{\nu_p}>0$.
The signed distance taken positive inside has gradient $-\nu_p$
at $p$. Differentiating along $\gamma$ at $0$ gives
\[
 d_\Om(\gamma(t))=at+o(t),\qquad
 |\gamma(t)-p|=O(t).
\]
Thus $\gamma$ is nontangential and the first part applies.
\end{proof}

\begin{remark}
Both distance estimates give squeezing convergence along the
entire inward normal. The remainder bound \eqref{eq:pointwise-taylor}
then gives strict pseudoconvexity by \Cref{cor:pointwise-taylor}.
With only local $C^2$ regularity, \Cref{thm:no-curves,thm:normal-log-growth}
instead rule out boundary curve germs and give the logarithmic
kernel growth exponent $n+1$ along the normal.
\end{remark}

\subsection{Biholomorphic realizations}
\label{subsec:realizations}

If a bounded $C^2$ strictly pseudoconvex domain is biholomorphic
to $\Om$, properness transfers its boundary squeezing limit to
$\Om$ without a boundary extension of the map. The local results
below apply to bounded targets; the unbounded example shows why
that restriction is needed.

\begin{cor}[Boundary restrictions on bounded realizations]
\label{cor:strongly-pseudoconvex-realization}
Let $D\Subset\C^n$, $n\ge2$, be strictly pseudoconvex with
$C^2$ boundary, and let $\Om\Subset\C^n$ be biholomorphic to
$D$. No boundary regularity of $\Om$ is needed for
$s_\Om(z)\to1$ as $d_\Om(z)\to0$.
\begin{enumerate}[label=\textup{(\roman*)},leftmargin=*]
\item If $\partial\Om$ is $C^2$ near $p\in\partial\Om$,
there is no boundary curve germ through $p$, and
\[
 \lim_{\delta\downarrow0}
 \frac{\log K_\Om(p-\delta\nu_p)}{\log(1/\delta)}=n+1.
\]
\item If $\partial\Om$ is $C^{2,\alpha}$ near $p$, for some
$0<\alpha\le1$, then $p$ is strictly pseudoconvex.
\end{enumerate}
In particular, the conclusions apply to bounded realizations
of the unit ball. If the entire boundary of $\Om$ is
$C^{2,\alpha}$, then $\Om$ is strictly pseudoconvex.
\end{cor}

\begin{proof}
First transfer the squeezing limit, without boundary extension.
Fix $\eta>0$. The set $\{w\in D:d_D(w)\ge\eta\}$ is compact
because $D$ is bounded. Its inverse image under the biholomorphism
$f:\Om\to D$ is compact in $\Om$ and, if nonempty, has positive
distance from $\partial\Om$. Thus $d_D(f(z))\to0$ as $d_\Om(z)\to0$.
Since $D$ has $C^2$ strictly pseudoconvex boundary,
\cite[Theorem~1.3]{DGZ} gives $s_D(f(z))\to1$. By biholomorphic
invariance, $s_\Om(z)=s_D(f(z))\to1$.

Biholomorphic invariance also makes $\Om$ pseudoconvex.
The transferred limit gives \eqref{eq:normal-hypothesis} in both cases:
use \Cref{thm:no-curves,thm:normal-log-growth} in (i)
and \Cref{thm:main} in (ii).
\end{proof}

Equivalently, a bounded pseudoconvex domain with a locally
$C^{2,\alpha}$ Levi-degenerate boundary point cannot be
biholomorphic to any bounded $C^2$ strictly pseudoconvex domain.

At a locally $C^2$ boundary point, a boundary curve germ or a
lower logarithmic growth exponent strictly below $n+1$ along
the normal rules out such a realization. In particular, the existence of the polydisc sequence in
\Cref{cor:polydisc-obstruction} rules out a biholomorphism to
any bounded $C^2$ strictly pseudoconvex domain.

\begin{remark}\label{rem:realization-priority}
The preceding realization argument uses properness rather than
boundary extension. For the earlier results of Pinchuk, Bell, Zimmer,
and Hu, see \Cref{subsec:earlier-results}. The necessary
conditions at $C^2$ regularity do not assert strict pseudoconvexity and do not
resolve the ball-image question in \cite[Problem~4.5]{FW}.
\end{remark}

Pinchuk \cite[p.~85]{Pin81}
gave an unbounded ball image with real-analytic boundary containing
a complex line. For $n=2$, the domain below is biholomorphically
equivalent to his example via $(z,w)\mapsto(z,-2w/(1+w))$.
In arbitrary dimension, its boundary is, after a rotation of $w$,
the positive-definite hyperquadric blow-up in
\cite[Example~2.4, p.~1706]{Kru20}; see also
\cite[Example~5.5]{KS13} for the two-dimensional model.
We include the explicit maps and kernel calculation to illustrate
the boundedness requirement in the realization result and the
normal kernel estimate.

\begin{prop}[A classical unbounded ball image]\label{prop:unbounded-ball-image}
For every $n\ge2$, the domain
\begin{equation}\label{eq:unbounded-ball-image}
 \Om_\infty=
 \{(z,w)\in\C^{n-1}\times\C:
       \Rea w>\norm{z}^2|w|^2\}
\end{equation}
is unbounded and biholomorphic to $\BB^n$. Its entire boundary
is a smooth real-analytic hypersurface, but its Levi form
vanishes on the complex tangent space at every point of
$M=\C^{n-1}\times\{0\}$.
Consequently, boundedness of the target in
\Cref{cor:strongly-pseudoconvex-realization} cannot be omitted,
even when the boundary is real-analytic everywhere.
\end{prop}

\begin{proof}
We check the boundary directly and obtain the ball map via
the Siegel model. Set $\rho(z,w)=\norm{z}^2|w|^2-\Rea w$.
If $z=0$ or $w=0$, then $\rho_w=-1/2$; otherwise
$\partial_z\rho=|w|^2\overline z\ne0$. The nonvanishing differential
makes $\rho$ change sign at every zero, so $\{\rho=0\}$ is exactly the
boundary of $\{\rho<0\}$ and is real-analytic.
At $w=0$, $\partial\rho=-dw/2$, so the complex tangent space
of the boundary is $T^{1,0}M$. Since $\rho|_M=0$, the tangential
Levi form vanishes. The ray
$\{(0,t):t>0\}$ proves unboundedness.
The domain inequality is the pullback of \eqref{eq:Siegel-domain}
by $P(z,w)=(wz,w)$, with holomorphic inverse $(Z,w)\mapsto(Z/w,w)$ since $w\ne0$ in the model.
Composing $P$ with \eqref{eq:inverse-cayley} and the unitary
change of sign in the last component gives
\begin{equation}\label{eq:unbounded-ball-map}
 \Psi(z,w)=
 \left(\frac{2wz}{1+w},\frac{1-w}{1+w}\right),
 \qquad
 \Psi^{-1}(\xi,\eta)=
 \left(\frac{\xi}{1-\eta},\frac{1-\eta}{1+\eta}\right).
\end{equation}
\end{proof}

The example shows where boundedness is used. The embedding definition
\eqref{eq:squeezing} applies to $\Om_\infty$ and gives
$s_{\Om_\infty}\equiv1$. For the defining polynomial $\rho(z,w)=\norm{z}^2|w|^2-\Rea w$
and for $\delta>0$, the automorphism
$A_\delta(z,w)=(\delta^{-1/2}z,\delta w)$ satisfies
$\rho\circ A_\delta=\delta\rho$. Hence
\begin{equation}\label{eq:unbounded-inverse-charts}
 G_\delta(\xi,\eta)
 :=A_\delta\circ\Psi^{-1}(\xi,\eta)
 =\left(\frac{\xi}{\sqrt\delta(1-\eta)},
         \delta\frac{1-\eta}{1+\eta}\right)
\end{equation}
is an inverse ball chart with $G_\delta(0)=(0,\delta)$.
For every fixed $\xi\ne0$ with $\norm{\xi}<1$,
$|G_\delta(\xi,0)|=|(\delta^{-1/2}\xi,\delta)|\to\infty$
as $\delta\downarrow0$. Thus the inverse charts
are not locally bounded, even though their centers approach
a finite real-analytic boundary point.

The Jacobian determinant of $P$ is $w^{n-1}$, so \eqref{eq:kernel-transform} and
the Siegel kernel \eqref{eq:Siegel-kernel} give
\begin{equation}\label{eq:unbounded-normal-kernel}
 K_{\Om_\infty}(z,w)
 =\frac{n!}{4\pi^n}
   \frac{|w|^{2n-2}}
        {(\Rea w-\norm{z}^2|w|^2)^{n+1}},
 \qquad
 K_{\Om_\infty}(0,\delta)
 =\frac{n!}{4\pi^n}\delta^{n-3}.
\end{equation}
The ball of radius $1/2$ centered at $(0,1/2)$ is contained
in $\Om_\infty$: its defining inequality implies
$\Rea w>|w|^2+|z|^2$ and $|w|<1$, hence
$\Rea w>|z|^2|w|^2$. At the origin the outward normal
is $-e_n$, by $d\rho=-d\Rea w$ there. The depth comparison
in \Cref{lem:normal-ball-comparison} does not require a bounded
ambient domain and gives
$d_{\Om_\infty}(0,\delta)=\delta$ for $0<\delta<1/2$.
The normal formula in \eqref{eq:unbounded-normal-kernel} gives
\[
 \delta^2K_{\Om_\infty}(0,\delta)
 =\frac{n!}{4\pi^n}\delta^{n-1}\longrightarrow0
 \qquad(n\ge2).
\]
Biholomorphism to a bounded domain therefore does not suffice
for the lower bound $c\delta^{-2}$ in \Cref{lem:kernel-bounds}.
The example is nonconvex: $(2e_1,1/8)$ and $(0,2)$ lie
in $\Om_\infty$, but their midpoint $(e_1,17/16)$ does not,
where $e_1$ is the first coordinate vector of $\C^{n-1}$.
It therefore does not conflict with Zimmer's theorem for
unbounded convex targets \cite[Theorem~1.1 and Remark~1.2]{Zim}.

\subsection{Examples and limits of the approach conditions}

The first example shows that the absence of boundary curves does
not by itself force the normal squeezing limit to be one. We then
contrast unrestricted tangential sequences with the normal
sampling criterion.

Put $\psi(0)=0$, $\psi(z)=\exp(-1/|z|^2)$ for $z\ne0$, and
consider
\[
 D=\{(z,w)\in\C^2: |z|<1/2,\ |w|<1,
                    \Rea w+\psi(z)<0\}.
\]
This is a bounded connected pseudoconvex domain: $\psi$ is
subharmonic on $|z|<1/2$, and the displayed defining
inequalities give a pseudoconvex intersection. Its boundary
has corners away from $p=(0,0)$, but is $C^\infty$ near $p$.
The point $p$ is Levi-degenerate and of infinite type.

No nonconstant boundary curve germ passes through $p$.
If $\gamma(\zeta)=(a(\zeta),b(\zeta))$ is such a germ, then
$\Rea b+\psi(a)=0$. On the set where $a\ne0$,
\[
 0=\partial_\zeta\partial_{\overline\zeta}(\psi\circ a)
   =\psi_{z\bar z}(a)|a'|^2,\qquad
 \psi_{z\bar z}(z)
   =e^{-1/|z|^2}\frac{1-|z|^2}{|z|^6}>0.
\]
A nonconstant $a$ would make the right-hand side positive
on a nonempty open set. Hence $a\equiv0$, and the open mapping
principle gives $b\equiv0$ as well.
Since $p$ is Levi-degenerate, the contrapositive of \Cref{thm:main} gives
\[
 \liminf_{\delta\downarrow0}s_D(0,-\delta)<1.
\]
The model is convex: for $0<r<1/2$ the radial and angular Hessian
eigenvalues of $\psi$ are $e^{-1/r^2}(4-6r^2)/r^6$ and
$2e^{-1/r^2}/r^4$, both positive. At the origin the Hessian is zero. The squeezing conclusion is
a lower-limit bound, not a uniform gap below $1$ for all small
$\delta$.

Unrestricted tangential sequences behave differently.
For example,
\[
 E=\{(z_1,z_2)\in\C^2:|z_1|^4+|z_2|^2<1\}
\]
has tangential sequences approaching its Levi-degenerate point $(0,1)$
along which squeezing tends to one \cite[Example~1.2, pp.~832--833]{NNN}.
Thus squeezing convergence on an unrestricted approaching sequence
does not imply strict pseudoconvexity.

For normal samples, \Cref{cor:normal-sampling} gives a direct
test: if $p$ is Levi-degenerate
and the boundary is $C^{2,\alpha}$ near $p$ for some
$0<\alpha\le1$, then for every fixed $0<\lambda<1$,
\[
 \liminf_{j\to\infty}
 s_\Om(p-\lambda^j\nu_p)<1.
\]
Otherwise, \Cref{cor:normal-sampling} would imply that $p$
is strictly pseudoconvex. The gap below $1$ in this lower limit
may depend on the domain, the point, and $\lambda$.
\section*{Acknowledgments}

The author thanks Professor Fusheng Deng for helpful discussions on related questions.

\section*{AI Use Disclosure}

The author conceived and developed the core ideas and proof strategy of this paper. ChatGPT was used as an aid in exploring possible constructions and checking the proofs. The author has independently verified the entire manuscript, including its mathematical arguments, references, and conclusions, and takes full responsibility for its contents.

\end{document}